\documentclass{amsart}
\usepackage{amssymb,cite,xfrac,fullpage,color}

\newcommand{\eq}[2]{\begin{equation}\label{eq:#1}#2\end{equation}}
\newcommand{\est}[2]{\begin{equation}\label{est:#1}#2\end{equation}}
\newcommand{\defn}[2]{\begin{equation}\label{defn:#1}#2\end{equation}}
\newcommand{\pde}[2]{\ensuremath{\left\{\begin{array}{l}#1\vspace{0.3cm}\\#2\end{array}\right.}}

\newcommand{\R}{\ensuremath{\mathbb{R}}}

\newcommand{\Z}{\ensuremath{\mathbb{Z}}}
\newcommand{\Schwartz}{\ensuremath{\mathcal{S}}}

\DeclareMathOperator{\supp}{supp}
\DeclareMathOperator{\sgn}{sgn}
\newcommand{\J}{L}
\newcommand{\X}{X}
\newcommand{\XEquiv}{\widetilde X}
\DeclareMathOperator{\Ai}{Ai}
\newcommand{\F}{\mathcal{F}}
\newcommand{\xiv}{\xi_v}
\renewcommand\Re{\operatorname{Re}}
\newcommand{\W}{\mathcal{W}}
\newcommand{\uasymp}{u_{\mathrm{asymp}}}
\newcommand{\uapp}{u_{\mathrm{app}}}
\newcommand{\DR}{\Omega^+}
\newcommand{\SSR}{\Omega^0}
\newcommand{\OR}{\Omega^-}
\newcommand{\err}{\mathbf{err}}
\newcommand{\VOR}{\mathbf\Omega^-}
\newcommand{\FOR}{\widehat\Omega^-}

\newcommand{\K}{Z}
\newcommand{\Y}{Y}
\newcommand{\udiag}{w_{v,+}}
\newcommand{\hyp}{\mathrm{hyp}}
\newcommand{\el}{\mathrm{ell}}
\newcommand{\In}{\mathrm{in}}
\newcommand{\Mid}{\mathrm{mid}}
\newcommand{\Out}{\mathrm{out}}

\newtheorem{lem}{Lemma}[section]
\newtheorem{thrm}[lem]{Theorem}
\newtheorem{propn}[lem]{Proposition}
\newtheorem{cor}[lem]{Corollary}
\theoremstyle{remark}
\newtheorem{rk}[lem]{Remark}
\newcommand{\lemma}[2]{\begin{lem}\label{lem:#1}#2\end{lem}}
\newcommand{\theorem}[2]{\begin{thrm}\label{thrm:#1}#2\end{thrm}}
\newcommand{\optheorem}[3]{\begin{thrm}[#2]\label{thrm:#1}#3\end{thrm}}
\newcommand{\proposition}[2]{\begin{propn}\label{propn:#1}#2\end{propn}}
\newcommand{\corollary}[2]{\begin{cor}\label{cor:#1}#2\end{cor}}
\newcommand{\remark}[2]{\begin{rk}\label{rk:#1}#2\end{rk}}
\newcommand{\bpf}{\begin{proof}}
\newcommand{\epf}{\end{proof}}

\title{Long time behavior of solutions to the mKdV}
\author[B. Harrop-Griffiths]{Benjamin Harrop-Griffiths}
\address{Department of Mathematics, University of California, Berkeley, CA 94720}
\email{benhg@math.berkeley.edu}
\thanks{The author was partially supported by the NSF grant DMS-1266182.}

\numberwithin{equation}{section}

\begin{document}

\maketitle


\begin{abstract}
In this paper we consider the long time behavior of solutions to the modified Korteweg-de Vries equation on \(\R\). For sufficiently small, smooth, decaying data we prove global existence and derive modified asymptotics without relying on complete integrability. We also consider the asymptotic completeness problem. Our result uses the method of testing by wave packets, developed in the work of Ifrim and Tataru on the \(1d\) cubic nonlinear Schr\"odinger and \(2d\) water wave equations.
\end{abstract}


\section{Introduction}

In this article we consider the long-time behavior of solutions to the mKdV equation
\eq{mKdV}{
\pde{
u_t+\tfrac{1}{3}u_{xxx}=\sigma (u^3)_x,\quad u\colon\R\times\R\rightarrow\R,
}{
u(0)=u_0,
}
}
where \(\sigma=\pm1\) and \(u_0\) is sufficiently small, smooth and decaying data.

The Cauchy problem for \eqref{eq:mKdV} has been studied extensively. For a summary of known results we refer the reader to \cite{LinaresPonce2009}. In particular, the mKdV is both locally well-posed \cite{Kato1983,KenigPonceVega1989,KenigPonceVega1993} and globally well-posed \cite{CollianderKeelStaffilaniTakaokaTao2003,Guo2009,Kishimoto2009} in \(H^s\) for \(s\geq\frac{1}{4}\). Below \(s=\frac{1}{4}\) the solution map fails to be uniformly continuous \cite{KenigPonceVega2001,ChristCollianderTao2003} although weaker forms of well-posedness hold \cite{ChristHolmerTataru2012}. Local well-posedness in non-\(L^2\)-based spaces closer to the critical scaling has also been obtained \cite{Grunrock2004,GrunrockVega2009}.

As the mKdV is completely integrable, global existence and asymptotic behavior can be studied using inverse scattering techniques such as in Deift and Zhou \cite{DeiftZhou1993} and references therein. A natural question to ask is whether it is possible to study the asymptotic behavior of the mKdV without relying on the completely integrable structure. Hayashi and Naumkin \cite{HayashiNaumkin1999,HayashiNaumkin2001} were able to prove global existence and derive modified asymptotics in a neighbourhood of a self-similar solution without relying on the complete integrability, with errors bounded in \(L^p\) for \(4<p\leq\infty\). Our result presents a significant improvement by proving modified scattering  in \(L^2\cap L^\infty\). We also derive the leading asymptotic in the oscillatory region and use slightly weaker assumptions on the initial data.

In the related case of the cubic nonlinear Schr\"odinger equation on \(\R\), modified asymptotics have been proved without inverse scattering techniques using both spatial methods \cite{LindbladSoffer2006} and  Fourier methods \cite{HayashiNaumkin1998NLS, KatoPusateri2011}. In this paper we use the \textit{method of testing by wave packets}, based on the work of Ifrim and Tataru on the \(1d\) cubic NLS \cite{IfrimTataru2014a} and \(2d\) water wave \cite{IfrimTataru2014b,IfrimTataru2014c} equations. This method essentially interpolates between the spatial side and Fourier side approaches by localizing in both space and frequency at the scale of the uncertainty principle.

In order to give a more complete picture of the asymptotic behavior of solutions, we consider the reciprocal problem: given a function with a suitable asymptotic profile, can we construct a solution to \eqref{eq:mKdV} matching this asymptotic behavior as \(t\rightarrow+\infty\)? Hayashi-Naumkin \cite{HayashiNaumkin2006} showed that under strong conditions on the data, including that it has mean zero, it is possible to find such a solution. Our result holds for a much larger class of data, including those with non-trivial mean. In the case of the gKdV, where solutions scatter to free solutions, asymptotic completeness was established by Cote \cite{Cote2006} and refined by Farah-Pastor \cite{FarahPastor2013}. Similar results have also been obtained for the cubic NLS, see for example \cite{IfrimTataru2014a} and references therein.

As in the case of the NLS \cite{DeiftZhou2002,HayashiNaumkin1998NLS,LindbladSoffer2006,IfrimTataru2014a}, Theorem \ref{thrm:Main} is also true for short-range perturbations of the form
\eq{PerturbedmKdV}{
\pde{
u_t+\tfrac{1}{3}u_{xxx}=(\sigma u^3+F(u))_x,
}{
u(0)=u_0,
}
}
where \(F\in C^2(\R)\) satisfies
\est{FEst}{
|F(u)|=O(|u|^p),\qquad |u|\rightarrow0,\quad p>3,
}
with some minor modifications if \(p\in(3,\tfrac{7}{2})\). For completeness we briefly outline these modifications in Appendix \ref{app:perturbations}.

While preparing this paper we learned that some similar results have been obtained by Germain, Pusateri and Rousset \cite{2015arXiv150309143G}.

For \(t>0\), the solution to the linear KdV equation
\eq{LinearKdV}{
\pde{
u_t+\tfrac{1}{3}u_{xxx}=0,
}{
u(0)=u_0,
}
}
is given by
\defn{FundamentalSolution}{
u(t,x)=t^{-\frac{1}{3}}\int\Ai(t^{-\frac{1}{3}}(x-y))u_0(y)\,dy,
}
where the Airy function is defined by the oscillatory integral
\[
\Ai(x)=\frac{1}{2\pi}\int e^{i(\frac{1}{3}\xi^3+x\xi)}\,d\xi.
\]

The linear KdV has Hamiltonian \(h(\xi)=-\frac{1}{3}\xi^3\) and hence the Hamiltonian flow associated to the linear KdV operator is
\eq{HFlo}{
(x,\xi)\mapsto(x-t\xi^2,\xi).
}
In particular, given a speed \(v\geq 0\), we expect wave packets with initial data localized in phase space near \((x,\xi)=(0,\pm\sqrt v)\) to travel along the ray \(\Gamma_v=\{x+tv=0\}\). As all wave packets travel towards \(x=-\infty\), we expect the solution to \eqref{eq:LinearKdV} to decay rapidly as \(t^{-\frac{1}{3}}x\rightarrow+\infty\) and oscillate as \(t^{-\frac{1}{3}}x\rightarrow-\infty\).

For Schwartz initial data, we can roughly divide the asymptotic behavior of solutions to \eqref{eq:LinearKdV} into three distinct regions as \(t\rightarrow+\infty\). In the decaying region \(t^{-\frac{1}{3}}x\rightarrow+\infty
\),
\[
u(t,x)=O(t^{-\frac{1}{3}}(t^{-\frac{1}{3}}|x|)^{-N}).
\]
In the self-similar region \(t^{-\frac{1}{3}}|x|\lesssim1\),
\[
u(t,x)=t^{-\frac{1}{3}}\Ai(t^{-\frac{1}{3}}x)\int u_0\,dx+O(t^{-\frac{2}{3}}).
\]
In the oscillatory region \(t^{-\frac{1}{3}}x\rightarrow-\infty\),
\[
u(t,x)=\pi^{-\frac{1}{2}}t^{-\frac{1}{3}}(t^{-\frac{1}{3}}|x|)^{-\frac{1}{4}}\Re\left(e^{i\phi}\hat u_0(t^{-\frac{1}{2}}|x|^{\frac{1}{2}})\right)+O(t^{-\frac{1}{3}}(t^{-\frac{1}{3}}|x|)^{-\frac{7}{4}}),
\]
where the phase is given by
\defn{phi}{
\phi(t,x)=-\frac{2}{3}t^{-\frac{1}{2}}|x|^{\frac{3}{2}}+\frac{\pi}{4}.
}

From \eqref{defn:FundamentalSolution}, we observe that if our initial data satisfies \(\|u_0\|_{{H^{0,1}}}\leq\epsilon\), the linear solution satisfies the dispersive estimates
\est{AiryEsts}{
|u(t,x)|\lesssim \epsilon t^{-\frac{1}{3}}\langle t^{-\frac{1}{3}}x\rangle^{-\frac{1}{4}},\qquad|u_x(t,x)|\lesssim \epsilon t^{-\frac{2}{3}}\langle t^{-\frac{1}{3}} x\rangle^{\frac{1}{4}},
}
and in particular,
\[
|uu_x|\lesssim \epsilon^2t^{-1}.
\]
We expect solutions to the nonlinear equation \eqref{eq:mKdV} to behave like solutions to the linear equation for sufficiently short times. So, if our initial data is of size \(\epsilon>0\) in a suitable norm, we expect it to satisfy \eqref{est:AiryEsts}, at least for sufficiently small \(T>1\). In particular, if \(\|\cdot\|\) is a Sobolev-type norm in \(x\) then
\[
\|u(t)\|\lesssim\|u(1)\|+\epsilon^2\int_1^t\|u(s)\|\,\frac{ds}{s}.
\]
The integral is bounded by \(\sup_t\|u(t)\|\) up to time \(T\approx e^{\epsilon^{-2}}\) and hence we only expect linear behavior up to this time. So while we may still have a global solution, we expect the asymptotic behavior of the solution to differ from the linear solution by a logarithmic difference in \(t\).

Our first result is that this is indeed the case.


\bigskip
\theorem{Main}{
There exists \(\epsilon>0\) such that for all \(u_0\in H^{1,1}\) satisfying
\est{InitialDataBnd}{
\|u_0\|_{H^{1,1}}\leq\epsilon,
}
there exists a unique global solution \(u\) to \eqref{eq:mKdV} with \(S(-t)u\in C(\R;H^{1,1})\) satisfying the estimates for \(t\geq1\) and a.e. \(x\in\R\),
\est{GlobalDecay}{
|u(t,x)|\lesssim \epsilon t^{-\frac{1}{3}}\langle t^{-\frac{1}{3}}x\rangle^{-\frac{1}{4}},\qquad |u_x(t,x)|\lesssim \epsilon {t^{-\frac{2}{3}}}\langle t^{-\frac{1}{3}}x\rangle^{\frac{1}{4}}.
}
Further, we have the following asymptotics as \(t\rightarrow+\infty\).

In the decaying region \(\DR_\rho=\{x>0:t^{-\frac{1}{3}}x\gtrsim t^{2\rho}\}\) we have the estimates
\est{RapidlyDecayingAsymptotics}{
\|t^{\frac{1}{3}}(t^{-\frac{1}{3}}x)^{\frac{3}{4}}u\|_{L^\infty(\DR_\rho)}\lesssim\epsilon,\qquad \|t^{\frac{1}{6}}(t^{-\frac{1}{3}}x)u\|_{L^2(\DR_\rho)}\lesssim\epsilon.
}

In the self-similar region \(\SSR_\rho=\{x\in\R:t^{-\frac{1}{3}}|x|\lesssim t^{2\rho}\}\), where \(0\leq\rho\leq\frac{1}{3}(\frac{1}{6}-C\epsilon^2)\), there exists a solution \(Q(y)\) to the Painlev\'e II equation
\eq{PLII}{
yQ-Q_{yy}+3\sigma Q^3=0,
}
satisfying
\est{QEst}{
|Q(y)|\lesssim\epsilon,
}
and we have the estimates
\est{SelfSimSoln}{
\|u-t^{-\frac{1}{3}}Q(t^{-\frac{1}{3}}x)\|_{L^\infty(\SSR_\rho)}\lesssim \epsilon t^{-\frac{1}{2}(\frac{5}{6}-C\epsilon^2)},\quad\|u-t^{-\frac{1}{3}}Q(t^{-\frac{1}{3}}x)\|_{L^2(\SSR_\rho)}\lesssim \epsilon t^{-\frac{2}{3}(\frac{5}{12}-C\epsilon^2)}.
}

In the oscillatory region \(\OR_\rho=\{x<0:\,t^{-\frac{1}{3}}|x|\gtrsim t^{2\rho}\}\), there exists a unique (complex-valued) function \(W\) {satisfying \(W(\xi)=\overline W(-\xi)\)} such that for \(C>0\) sufficiently large,
\est{AEst}{
\|W\|_{H^{1-C\epsilon^2,1}\cap L^\infty}\lesssim\epsilon,
}
and
\est{QuasilinearAsymptotics}{
u(t,x)=\pi^{-\frac{1}{2}}t^{-\frac{1}{3}}(t^{-\frac{1}{3}}|x|)^{-\frac{1}{4}}\Re\left(e^{i\phi(t,x)+\frac{3 i\sigma}{4\pi}|W(t^{-\frac{1}{2}}|x|^{\frac{1}{2}})|^2\log({t^{-\frac12}|x|^{\frac32}})}W(t^{-\frac{1}{2}}|x|^{\frac{1}{2}})\right)+\err_x,
}
where the error satisfies the estimates
\est{ErrorEstSpace}{
\|t^{\frac13}(t^{-\frac{1}{3}}|x|)^{\frac{3}{8}}\err_x\|_{L^\infty(\OR_\rho)}\lesssim\epsilon,\qquad \|t^{\frac16}(t^{-\frac{1}{3}}|x|)^{\frac{1}{4}}\err_x\|_{L^2(\OR_\rho)}\lesssim\epsilon.
}

In the corresponding frequency region \(\FOR_\rho=\{\xi>0:t^{\frac{1}{3}}\xi\gtrsim t^\rho\}\) we have
\est{QuasilinearAsymptoticsFreq}{
\hat u(t,\xi)=e^{\frac{1}{3}it\xi^3+\frac{3i\sigma}{4\pi}|W(\xi)|^2\log({t\xi^3})}W(\xi)+\err_\xi,
}
where the error satisfies
\est{ErrorEstFreq}{
\|(t^{\frac{1}{3}}\xi)^{\frac{1}{4}}\err_\xi\|_{L^\infty(\FOR_\rho)}\lesssim\epsilon,\qquad\|t^{\frac{1}{6}}(t^{\frac{1}{3}}\xi)^{\frac{1}{2}}\err_\xi\|_{L^2(\FOR_\rho)}\lesssim\epsilon.
}
}
\bigskip

\remark{MainThrmNotes}{
As \eqref{eq:mKdV} has time reversal symmetry given by
\[
u(t,x)\mapsto u(-t,-x),
\]
we get corresponding asymptotics as \(t\rightarrow-\infty\).
}

\medskip

\remark{LossofReg}{
The loss of regularity of \(W\) in Theorem \ref{thrm:Main} can be compared to the similar results \cite{IfrimTataru2014a,IfrimTataru2014b}. Indeed, as the direct scattering problem for the cubic NLS and mKdV is the same, we expect the correspondence between the \(W\) of Theorem \ref{thrm:Main} and \(u_0\) to be the same as in Theorem 1 of \cite{IfrimTataru2014a}. From the inverse scattering theory, see for example \cite{DeiftZhou1993,DeiftZhou2002}, we expect this loss of regularity to be logarithmic in nature.
}
\medskip


For the asymptotic completeness, a key object of study will be the one-parameter family of solutions to the Painlev\'e II equation \eqref{eq:PLII}. We first state the following result giving the asymptotic behavior of these solutions.

\bigskip
\optheorem{PLII}{Deift-Zhou \cite{DeiftZhou1995}}{
Given \(W\in\R\) {(and sufficiently small if \(\sigma=-1\))} there exists a unique solution \(Q(y;W)\) to the Painlev\'e II equation \eqref{eq:PLII} such that
\est{PLIIAsymp-}{
\begin{gathered}
Q(y;W)=\pi^{-\frac{1}{2}}|y|^{-\frac{1}{4}}\Re\left(e^{-\frac23i|y|^{\frac32}+i\frac\pi4+\frac{3i\sigma}{4\pi}W^2\log|y|^{\frac{3}{2}}+i\sigma \theta(W^2)}W\right)+O(|y|^{-\frac{5}{4}}\log|y|),\quad y\rightarrow-\infty,
\\
Q(y;W)=q_\sigma(W)\Ai(y)+O(|y|^{-\frac{1}{4}}e^{-\frac{4}{3}y^{\frac{3}{2}}}),\quad y\rightarrow+\infty,
\end{gathered}
}
where
\begin{gather*}
\theta(W^2)=\tfrac{9\log 2}{4\pi}W^2-\arg\Gamma\left(\tfrac{3i}{4\pi}W^2\right)-\tfrac{\pi}{2},\qquad q_\sigma(W)=\sgn W\left(\tfrac{2\sigma}{3}\left(1-e^{-\frac{3\sigma}{2}W^2}\right)\right)^{\frac{1}{2}}.
\end{gather*}
}
\bigskip


For a real-valued {even function \(W\)}, we define
\defn{uasymp}{
\uasymp(t,x)=t^{-\frac{1}{3}}Q(t^{-\frac{1}{3}}x;W(t^{-\frac{1}{2}}|x|^{\frac{1}{2}})).
}
We observe that from Theorem \ref{thrm:PLII}, this has an asymptotic profile matching that of Theorem \ref{thrm:Main}. We then look for a solution to the problem
\eq{ACPDE}{
\pde{
u_t+\tfrac{1}{3}u_{xxx}=\sigma(u^3)_x,
}{
\lim\limits_{t\rightarrow+\infty}\left(u(t)-\uasymp(t)\right)=0.
}
}
We define the space \(\Y\) of real-valued even functions with norm
\defn{YNorm}{
\|W\|_\Y=\|\langle D\rangle^{C\epsilon^2}W\|_{H^{1,1}}
}
and then have the following asymptotic completeness result.

\bigskip
\theorem{AC}{
There exist \(\epsilon,C>0\) such that for all \(W\in\Y\) satisfying
\est{WAC}{
\|W\|_\Y\leq\epsilon,
}
there exists a unique solution to \eqref{eq:ACPDE} such that \(S(-t)u\in C(\R;H^{1,1})\).
}
\bigskip


\remark{LogarithmicDifference}{
Similar to Theorem \ref{thrm:Main} we have a loss of regularity between \(W\) and \(u\). In order to close the argument we require an extra \(C\epsilon^2\) derivatives for both \(W_z\) and \(zW\).
}

\medskip

\remark{SmallGap}{
As we use the 1-parameter family of real-valued solutions to the Painlev\'e II as our asymptotic object, we are restricted to considering real-valued W. This leaves a small gap between Theorems \ref{thrm:Main} and \ref{thrm:AC}.
}

\medskip


We conclude this section by giving an outline of the proof of Theorems \ref{thrm:Main} and \ref{thrm:AC}. In order to control the spatial localization of solutions we look to control the ``vector field"
\defn{GalileanVF}{
\J u=S(t)xS(-t)u=(x-t\partial_x^2)u.
}
However, \(\J \) does not behave well with respect to the nonlinearity, so as in \cite{HayashiNaumkin1998gKdV,HayashiNaumkin1999,HayashiNaumkin2001} we instead work with
\defn{ScalingVF}{
\Lambda u=\partial_x^{-1}(3t\partial_t+x\partial_x+1)u.
}
We observe that if \(u\) is a solution to \eqref{eq:mKdV} then
\eq{LambdaIsNF}{
\Lambda u=\J u+3t\sigma u^3.
}
As \(3t\partial_t+x\partial_x+1\) generates the mKdV scaling symmetry
\eq{Scaling}{
u(t,x)\mapsto\lambda u(\lambda^3t,\lambda x),\qquad u_0(x)\mapsto\lambda u_0(\lambda x),
}
the function \(v=\Lambda u\) satisfies the linearized equation
\eq{DiffByScaling}{
\pde{
v_t+\tfrac{1}{3}v_{xxx}=3\sigma u^2v_x,
}{
v(0)=xu_0.
}
}

For a large fixed constant \(M_0\geq2\) we define the space \(\X\) with norm
\defn{Xnorm}{
\|u\|_{\X}^2=\|u\|_{H^1}^2+\langle t\rangle^{-2\delta}\|\Lambda u\|_{L^2}^2,
}
where
\defn{delta}{
\delta=3M_0^2\epsilon^2.
}
We then have the following local well-posedness result that can be proved as in Kenig-Ponce-Vega \cite{KenigPonceVega1989,KenigPonceVega1993}.

\bigskip
\theorem{LocalWellPosedness}{
If \(u_0\in H^{1,1}\) satisfies \eqref{est:InitialDataBnd} then there exists \(T=T(\epsilon)\rightarrow\infty\) as \(\epsilon\rightarrow0\) and a unique solution \(u\in C([0,T];\X)\) such that
\est{LWPBnd}{
\sup\limits_{t\in[0,T]}\|u(t)\|_{\X}\leq 10\epsilon.
}
Further, the solution map \(u_0\mapsto u(t)\) is locally Lipschitz.
}
\bigskip

In \S\ref{sect:scattering} we prove Theorem \ref{thrm:Main}. Using the local well-posedness result, for \(\epsilon>0\) sufficiently small we can find \(T>1\) and a unique solution \(u\in C([0,T]; X)\) to \eqref{eq:mKdV}. We then make the bootstrap assumption that \(u\) satisfies the linear pointwise estimate
\est{Bootstrap}{
\sup\limits_{t\in[1,T]}\left(\|t^{\frac{1}{3}}\langle t^{-\frac{1}{3}}x\rangle^{\frac{1}{4}}u\|_{L^\infty}+\|t^{\frac{2}{3}}\langle t^{-\frac{1}{3}}x\rangle^{-\frac{1}{4}}u_x\|_{L^\infty}\right)\leq M_0\epsilon
}
and show that under this assumption, for \(\epsilon>0\) sufficiently small, we have the energy estimate
\est{GlobalNRG}{
\sup\limits_{t\in[0,T]}\|u\|_{\X}\lesssim\epsilon,
}
with a constant independent of \(M_0,T\). To complete the proof of global existence we need to close the bootstrap estimate \eqref{est:Bootstrap}.

To control the pointwise behavior of solutions we use the method of testing by wave packets \cite{IfrimTataru2014a,IfrimTataru2014b,IfrimTataru2014c}. A wave packet is an approximate solution localized in both space and frequency on the scale of the uncertainty principle. We define a wave packet \(\Psi_v\) adapted to the ray \(\Gamma_v\) and measure \(u\) along \(\Gamma_v\) by considering
\defn{gammaBasic}{
\gamma(t,v)=\int u(t,x)\overline\Psi_v(t,x)\,dx.
}

A key innovation of Ifrim and Tataru is to choose the wave packet to be localized at a \(t\)-dependent scale. For the KdV, a wave packet adapted to the ray \(\Gamma_v\) will be localized at scale \(\lambda=t^{-\frac13}\langle t^{\frac23}v\rangle^{-\frac14}\) in frequency and at scale \(\lambda^{-1}\) in space. However, as we only make use of the wave packets in the region
\defn{VOR}{
\VOR_0=\{v>0:t^{\frac{2}{3}}v\gtrsim 1\}
}
corresponding to the region \(\OR_0\), we instead define
\defn{lambda}{
\lambda=t^{-\frac{1}{2}}v^{-\frac{1}{4}}.
}
We then reduce closing the bootstrap estimate \eqref{est:Bootstrap} to proving global bounds for \(\gamma\). To derive these bounds, we show that \(\gamma\) satisfies an ODE of the form
\[
\dot\gamma(t,v)=3i\sigma t^{-1}|\gamma(t,v)|^2\gamma(t,v)+\textrm{error}.
\]
The logarithmic correction to the phase then arises as a consequence of solving this ODE.

In \S\ref{sect:AsympComplete} we prove Theorem \ref{thrm:AC}. The key idea here is to replace \(\uasymp\) by a regularized version \(\uapp\), where the regularization is on the scale of the wave packets. The result then follows by applying a contraction mapping argument to the resulting equation for the difference \(v=u-\uapp\) in a suitable space.


\section{Notation and definitions}\label{sect:note}

We recall that solutions to \eqref{eq:mKdV} have conserved quantities
\begin{gather}
E_0(t)=\int u\,dx,\label{est:ConservedMean}\\
E_1(t)=\int u^2\,dx,\label{est:ConservedMass}\\
E_2(t)=\int u_x^2+\tfrac{3}{2}\sigma u^4\,dx.\label{est:ConservedEnergy}
\end{gather}
We note that as \eqref{eq:mKdV} is completely integrable there are an infinite number of higher order conserved quantities.

We define the Fourier transform of a Schwartz function \(f\in\Schwartz(\R)\) to be
\defn{FT}{
\hat f(\xi)=\F f(\xi)=\int f(x)e^{-ix\xi}\,dx,
}
with inverse
\defn{IFT}{
\check f(x)=\F^{-1}f(x)=\frac{1}{2\pi}\int f(\xi)e^{ix\xi}\,d\xi.
}
The linear KdV propagator \(S(t)\) can then be written as
\defn{LinearKdVPropagator}{
S(t)f=\frac{1}{2\pi}\int\hat f(\xi)e^{i(\frac{1}{3}t\xi^3+x\xi)}\,d\xi.
}
{
We then have the dispersive estimate,
\est{Dispersive}{
\|S(t)f\|_{L^\infty}\lesssim t^{-\frac13}\|f\|_{L^1}.
}

Let \(\psi\in C^\infty_0\) be a real-valued, even function satisfying \(0\leq\psi\leq1\), supported on \((-2,2)\) and identically \(1\) on \([-1,1]\). For \(N\in 2^\Z\), we define \(\psi_N(\xi)=\psi(N^{-1}\xi)\) and \(\varphi_N(\xi)=\psi_N(\xi)-\psi_{\frac N2}(\xi)\). We then have the Littlewood-Paley projections
\defn{LWPProj}{
P_Nu=\varphi_N(D)u,\qquad P_{\leq N}u=\psi_N(D)u,\qquad P_{>N}u=u-P_{\leq N}u,\qquad P_{N<\cdot\leq M}u=P_{\leq M}P_{>N}u.
}
}
We also define the projections to positive and negative frequencies
\defn{CauchyOps}{
P_\pm u=\mathbf{1}_{(0,\infty)}(\pm D)u.
}
We recall the Bernstein inequality, for \(1\leq p\leq q\leq \infty\),
\est{Bernie}{
\|P_Nu\|_{L^q}\lesssim N^{\frac{1}{p}-\frac{1}{q}}\|P_Nu\|_{L^p},
}
and the Sobolev estimate
\est{SobolevEst}{
\|u\|_{L^\infty}\lesssim\|u\|_{L^2}^{\frac12}\|u_x\|_{L^2}^{\frac12}.
}

{
We recall that if \(\chi\in C^\infty_0\) and for \(R>0\) we define \(\chi_R(x)=\chi(R^{-1}x)\) then for \(1\leq p\leq\infty\) we have the estimate
\est{LocalisationHelper}{
\|(1-P_{\frac N4\leq\cdot\leq 4N})(\chi_RP_Nu)\|_{L^p}\lesssim_k \langle RN\rangle^{-k}\|P_Nu\|_{L^p},
}
so \(\chi_RP_Nu\) is localized at frequencies \(\sim N\) whenever \(RN\gg1\).
}
As a consequence we may obtain a version of Bernstein's inequality whenever \(RN\gg1\) and \(1\leq p\leq q\leq\infty\)
\begin{align*}
\|\chi_RP_Nu\|_{L^q} &\leq \|P_{\frac N4\leq\cdot\leq 4N}(\chi_RP_Nu)\|_{L^q} + \|(1-P_{\frac N4\leq\cdot\leq 4N})(\chi_RP_Nu)\|_{L^q}\\
&\lesssim_k N^{\frac1p - \frac1q}\|P_{\frac N4\leq\cdot\leq 4N}(\chi_RP_Nu)\|_{L^p} + \langle RN\rangle^{-k}\|P_N\|_{L^q}\\
&\lesssim_k N^{\frac1p - \frac1q}\left(\|\chi_RP_Nu\|_{L^p} + \langle RN\rangle^{-k}\|P_N\|_{L^p}\right)
\end{align*}

We define the weighted Sobolev norms \(H^{k,j}\) by
\defn{Hkj}{
\|u\|_{H^{k,j}}^2=\|\langle\xi\rangle^k\hat u\|_{L^2}^2+\|\langle x\rangle^j u\|_{L^2}^2.
}

We call a pair of indices \((p,q)\) admissible if
\defn{AdmissiblePair}{
\frac{2}{p}+\frac{1}{q}=\frac{1}{2},\quad 4\leq p\leq\infty,\quad 2\leq q\leq\infty.
}
If we define
\defn{BigPhi}{
\Phi f=\int_t^\infty S(t-s)f(s)\,ds,
}
then for admissible pairs \((p_1,q_1),(p_2,q_2)\) and any decomposition \(f=f_1+f_2\) we have the estimate \cite{KenigPonceVega1989,KenigPonceVega1993}
\est{KPV}{
\|\Phi f\|_{L^\infty_tL^2_x}+\||D|^{1-\frac{5}{{p_1}}}\Phi f\|_{L^{p_1}_xL^{q_1}_t}\lesssim\|f_1\|_{L^1_tL^2_x}+\||D|^{\frac{5}{p_2}-1}f_2\|_{L^{p_2'}_xL^{q_2'}_t}.
}
For \((p,q)\) admissible we have the estimate \cite[Corollary 3.6]{ChristHolmerTataru2012}
\est{UpEst}{
\|u\|_{L^\infty_tL^2_x}+\||D|^{1-\frac{5}{p}}u\|_{L^p_xL^q_t}\lesssim\|u\|_{U^{{\min\{p,q\}}}_S},
}
where the space \(U^p_S\) is defined as in \cite{KochTataruVisan2014}. For {\(p>2\) we have the embedding,
\est{UpVpEmbed}{
\dot W^{1,1}_tL^2_x\subset V^2_{rc}\subset U^p,
}
}
where \(V^p_{rc}\) is the space of right-continuous functions of bounded \(p\)-variation  (see for example \cite{KochTataruVisan2014}).

\section{Modified Scattering}\label{sect:scattering}


\subsection{Energy estimates}

We first derive energy estimates for \(u\) under the the bootstrap assumption \eqref{est:Bootstrap}. Our argument is similar to Hayashi-Naumkin \cite{HayashiNaumkin1998gKdV,HayashiNaumkin1999,HayashiNaumkin2001}.

\medskip
\proposition{NRG}{
For \(\epsilon>0\) chosen sufficiently small and \(t\in [0,T]\) we have the energy estimates
\est{SobolevNRG}{
\|u\|_{H^1}\lesssim\epsilon,
}
\est{JNRG}{
\|\Lambda u\|_{L^2}\lesssim\epsilon\langle t\rangle^\delta,
}
where \(\delta\) is defined as in \eqref{defn:delta} and the constants are independent of \(M_0,T\).
}
\bpf
From conservation of mass \eqref{est:ConservedMass}, we have
\[
\|u\|_{L^2}=\|u_0\|_{L^2}\leq\epsilon.
\]

From the Sobolev estimate \eqref{est:SobolevEst}, for any \(\theta>0\), we have
\[
\int u^4\,dx\lesssim \|u\|_{L^2}^3\|u_x\|_{L^2}\lesssim \theta^{-1}\|u\|_{L^2}^6+\theta\|u_x\|_{L^2}^2\lesssim \theta^{-1}\epsilon^4\|u\|_{L^2}^2+\theta\|u_x\|_{L^2}^2.
\]
Defining \(E_j(t)\) for \(j=1,2\) as in \eqref{est:ConservedMass} and \eqref{est:ConservedEnergy}, for \(\theta>0\) chosen sufficiently small we then have
\est{H1IsNRG}{
\|u(t)\|_{H^1}\sim E_1(t)+E_2(t)=E_1(0)+E_2(0)\sim\|u_0\|_{H^1}\leq\epsilon,
}
where the constants are independent of \(M_0\).

If \(v=\Lambda u\), from the estimate \eqref{est:LWPBnd} we have
\[
\sup\limits_{t\in[0,1]}\|v(t)\|_{L^2_x}\lesssim\epsilon.
\]
For \(t\geq1\), as a consequence of \eqref{est:Bootstrap} we have the estimate
\eq{uuxBnd}{
\|uu_x\|_{L^\infty}\leq M_0^2\epsilon^2t^{-1},\quad t\geq1,
}
and we may then use the equation \eqref{eq:DiffByScaling} to get
\begin{align*}
\partial_t\|v\|_{L^2}^2&=6\sigma\int u^2v_xv\,dx\\
&=-6\sigma\int uu_xv^2\,dx\\
&\leq 6M_0^2\epsilon^2{t^{-1}}\|v\|_{L^2}^2.
\end{align*}
Using Gronwall's inequality for \(t\geq1\) we obtain \eqref{est:JNRG}.

\epf
\medskip

For convenience we will make use of the norm
\[
\|u\|_{\XEquiv}^2=\|Lu\|_{L^2}^2+\|t^{\frac13}\langle t^{\frac13}D_x\rangle^{-1}u\|_{L^2}^2.
\]
For \(N\in 2^\Z\), we define
\[
u_N=P_Nu,\qquad u_{\leq t^{-\frac13}}=\sum\limits_{N\leq t^{-\frac13}}u_N.
\]
We note that we have the compatibility estimate
\est{Compat}{
\|u\|_{\XEquiv}^2\sim\|u_{\leq t^{-\frac13}}\|_{\XEquiv}^2+\sum\limits_{N>t^{-\frac13}}\|u_N\|_{\XEquiv}^2.
}
We then have the following corollary to Proposition \ref{propn:NRG}.

\medskip
\corollary{AltNorm}{
For \(\epsilon>0\) sufficiently small and \(t\in[1,T]\), we have the estimate
\est{AltNorm}{
\|u\|_{\XEquiv}\lesssim\epsilon t^{\frac16}.
}
}
\bpf
From the estimate \eqref{est:LWPBnd} and the bootstrap assumption \eqref{est:Bootstrap} we obtain
\est{uLpBnd}{
\|u\|_{L^p}\lesssim M_0\epsilon \langle t\rangle^{\frac{1}{3p}-\frac{1}{3}},\quad p\in(4,\infty].
}
From the energy estimate \eqref{est:JNRG} we then have
\[
\|\J u\|_{L^2}\lesssim \|\Lambda u\|_{L^2}+3t\|u\|_{L^6}^3\lesssim \epsilon \langle t\rangle^\delta+(M_0\epsilon)^3\langle t\rangle^{\frac16}.
\]
This gives the first component of \eqref{est:AltNorm}, provided \(\epsilon=\epsilon(M_0)>0\) is chosen sufficiently small that \(\delta\in(0,\frac{1}{6}]\).

For the second part we make a self-similar change of variables by defining
\defn{BigU}{
U(t,y)=t^{\frac13}u(t,t^{\frac13}y).
}
We observe that \(U\) satisfies the equation
\eq{mKdVForBigU}{
\pde{
\partial_tU=\frac13t^{-1}\partial_y(yU-U_{yy}+3\sigma U^3)
}{
U(1,y)=u(1,y)
}
}
From the energy estimate \eqref{est:JNRG} we have
\begin{align*}
\partial_t\|\langle D_y\rangle^{-1}U\|_{L^2_y} &\lesssim t^{-1}\|yU-U_{yy}+3\sigma U^3\|_{L^2_y}\\
&\lesssim t^{-\frac76}\|\Lambda u\|_{L^2}\\
&\lesssim\epsilon t^{\delta-\frac76}.
\end{align*}
So for \(\epsilon>0\) chosen sufficiently small that \(\delta<\frac16\) we may integrate to obtain
\[
\|\langle D_y\rangle^{-1}U\|_{L^2_y}\lesssim \epsilon,
\]
which gives us the second component of \eqref{est:AltNorm}.
\epf
\medskip

{
\subsection{Initial bounds}\label{sect:InitBnds}

In this section we prove a number of estimates for \(u\) that will allow us to reduce closing the bootstrap estimate \eqref{est:Bootstrap} to considering the behavior of \(u\) along the rays \(\Gamma_v\) for \(v\in\VOR_0\). Our argument uses Klainerman-Sobolev type estimates similar to \cite{HunterIfrimTataru2014,IfrimTataru2014b,IfrimTataru2014c}.

Let \(t\geq1\) be fixed. We first decompose \(u\) into a piece on which \(L\) acts hyperbolically and piece on which it acts elliptically. Let \(\psi\in C^\infty_0\) be defined as in \S\ref{sect:note}. Let \(\nu\gg1\) be a fixed parameter and define
\[
\chi(x)=\psi(\nu^{-1}x)-\psi(\nu x),\qquad \chi^\hyp=\mathbf1_{(-\infty,0)}\chi,\qquad \chi^\el=1-\chi^\hyp.
\]
We then rescale, defining {\(\chi_N(t,x)=\chi(t^{-1}N^{-2}x)\)} and similarly for \(\chi_N^\hyp\), \(\chi_N^\el\).

For each \(N>t^{-\frac13}\), we decompose \(u_N\) as
\[
u_N=u_{N,+}^\hyp+u_{N,-}^\hyp+u_N^\el,
\]
where \(u_{N,\pm}^\hyp=\chi_N^\hyp P_\pm u_N\).
We then define the hyperbolic parts of \(u\) by
\[
u^\hyp_\pm=\sum\limits_{N> t^{-\frac13}}u_{N,\pm}^\hyp,
\]
and use this to decompose \(u\),
\[
u=u^\hyp_++u^\hyp_-+u^\el.
\]
We observe that \(u^\hyp=u^\hyp_++u^\hyp_-=2\Re(u^\hyp_+)\) is supported in the region \(\{t^{-\frac13}x<-\nu^{-1}\}\) and we may define \(\OR_0\) such that \(\OR_0\subset\{t^{-\frac13}x<-\nu^{-1}\}\).

In the region \(\OR_0\), the symbol of \(L\) factorizes as
\[
x-t\xi^2=-(|x|^{\frac12}\mp t^{\frac12}\xi)(|x|^{\frac12}\pm t^{\frac12}\xi),
\]
and hence we define operators associated to this factorization,
\[
L_\pm=|x|^{\frac12}\pm it^{\frac12}\partial_x.
\]
We note that \(L_-\) is elliptic on positive frequencies and \(L_+\) is elliptic on negative frequencies.

We then have the following bounds for the hyperbolic and elliptic parts of \(u\).

\medskip
\proposition{InitialPointwiseBounds}{
For \(t\in[1,T]\) we may decompose \(u=u^\hyp+u^\el\) into a hyperbolic part \(u^\hyp\) supported in \(\OR_0\) and an elliptic part \(u^\el\), satisfying the bounds
\begin{gather}
\|t^{\frac16}\langle t^{-\frac13}x\rangle u^\el\|_{L^2}\lesssim\epsilon,\qquad\|t^{\frac13}\langle t^{-\frac13}x\rangle^{\frac34}u^\el\|_{L^\infty}\lesssim\epsilon,\qquad \|t^{\frac23}\langle t^{-\frac13}x\rangle^{\frac14}u^\el_x\|_{L^\infty}\lesssim\epsilon,\label{est:EllipticBounds}\\
\|t^{\frac13}u^\hyp\|_{L^\infty}\lesssim\epsilon,\qquad \|t^{\frac23}\langle t^{-\frac13}x\rangle^{-\frac12}u^\hyp_x\|_{L^\infty}\lesssim\epsilon\label{est:HyperbolicBounds}.
\end{gather}
}
\medskip

In order to prove Proposition \ref{propn:InitialPointwiseBounds}, we first prove the following elliptic estimates for our solution.

\medskip
\lemma{EllipticEstimates}{
For \(t\in[1,T]\) we have the estimates
\begin{align}
\|t^{\frac13}\langle t^{-\frac13}x\rangle u_{\leq t^{-\frac13}}\|_{L^2} &\lesssim \|u_{\leq t^{-\frac13}}\|_{\XEquiv},\label{est:Low}\\
\|(|x|+tN^2)u_N^\el\|_{L^2} &\lesssim\|u_N\|_{\XEquiv}\label{est:Elliptic},\qquad N> t^{-\frac13},\\
\|(|x|^{\frac12}+t^{\frac12}N)L_\pm u_{N,\pm}^\hyp\|_{L^2} &\lesssim\|u_N\|_{\XEquiv}\label{est:Hyperbolic},\qquad N> t^{-\frac13}.
\end{align}
}
\bpf~

\textit{A. Low frequencies.} Using the frequency localization of \(u_{\leq t^{-\frac13}}\) wo obtain
\begin{gather*}
\|t^{\frac13}u_{\leq t^{-\frac13}}\|_{L^2}\lesssim \|t^{\frac13}\langle t^{\frac13}D_x\rangle^{-1}u_{\leq t^{-\frac13}}\|_{L^2},\\
\|xu_{\leq t^{-\frac13}}\|_{L^2}\lesssim \|Lu_{\leq t^{-\frac13}}\|_{L^2}+\|t^{\frac13}\langle t^{\frac13}D_x\rangle^{-1}u_{\leq t^{-\frac13}}\|_{L^2},
\end{gather*}
which give us the estimate \eqref{est:Low}.

\textit{B. Elliptic region.} Let \(N> t^{-\frac13}\). By rescaling under the mKdV scaling \eqref{eq:Scaling}, it suffices to consider the case \(N=1\). We decompose
\begin{gather*}
\chi_1^\el=\chi_1^\In+\chi_1^\Out+\chi_1^\Mid,\\
\chi_1^\In(t,x)=\psi(\nu t^{-1}x),\qquad \chi_1^\Out(t,x)=1-\psi(\nu^{-1} t^{-1}x),\qquad\chi_1^\Mid(t,x)=\chi_1^\hyp(t,-x).
\end{gather*}
The functions \(\chi_1^\In,\chi_1^\Mid\in C^\infty_0\) and \(\chi_1^\Out\in C^\infty\) are supported in the regions \(\{|x|<2\nu^{-1}t\}\), \(\{\nu^{-1}t<x<2\nu t\}\) and \(\{|x|>\nu t\}\) respectively.

\textit{B(i). Inner region.}
Using the estimate \eqref{est:LocalisationHelper}, we have
\begin{align*}
t\|\chi_1^\In u_1\|_{L^2}&\lesssim t\|P_{\frac 14\leq\cdot\leq 4}(\chi_1^\In u_1)\|_{L^2}+t\|(1-P_{\frac 14\leq\cdot\leq 4})(\chi_1^\In u_1)\|_{L^2}\\
&\lesssim_k t\|\partial_x^2P_{\frac 14\leq\cdot\leq 4}(\chi_1^\In u_1)\|_{L^2}+t \langle \nu^{-1}t\rangle^{-k}\|u_1\|_{L^2}\\
&\lesssim t\|\chi_1^\In \partial_x^2u_1\|_{L^2}+C(\nu)\|u_1\|_{\XEquiv}.
\end{align*}
As a consequence we obtain
\begin{align*}
t\|\chi_1^\In u_1\|_{L^2}&\lesssim \|(x-t\partial_x^2)u_1\|_{L^2}+C(\nu)\|u_1\|_{\XEquiv}+\|x\chi_1^\In u_1\|_{L^2}\\
&\lesssim \|(x-t\partial_x^2)u_1\|_{L^2}+C(\nu)\|u_1\|_{\XEquiv}+\nu^{-1}t\|\chi_1^\In u_1\|_{L^2}.
\end{align*}
Taking \(\nu\gg1\) to be a sufficiently large fixed constant, we have
\[
\|(|x|+t)\chi_1^\In u_1\|_{L^2}\lesssim t\|\chi_1^\In u_1\|_{L^2}\lesssim \|u\|_{\XEquiv},
\]
which gives us \eqref{est:Elliptic} in the region \(\{|x|<2\nu^{-1}t\}\).

\textit{B(ii). Outer region.}
Proceeding as for the inner region we have
\begin{align*}
\|x\chi_1^\Out u_1\|_{L^2}&\lesssim \|(x-t\partial_x^2)u_1\|_{L^2}+C(\nu)\|u_1\|_{\XEquiv}+t\|\chi_1^\Out u_1\|_{L^2}\\
&\lesssim \|(x-t\partial_x^2)u_1\|_{L^2}+C(\nu)\|u_1\|_{\XEquiv}+\nu^{-1}\|x\chi_1^\Out u_1\|_{L^2},
\end{align*}
and hence
\[
\| (|x|+t)\chi_1^\Out u_1\|_{L^2}\lesssim \|x\chi_1^\Out u_1\|_{L^2}\lesssim \|u\|_{\XEquiv},
\]
which gives us \eqref{est:Elliptic} in the region \(\{|x|>\nu t\}\).

\textit{B(iii). Middle region.} We now ignore the dependence of constants upon \(\nu\). Integrating by parts we have
\begin{align*}
&\|\chi_1^\Mid x u_1\|_{L^2}^2+\|\chi_1^\Mid t\partial_x^2u_1\|_{L^2}^2+2t\int (\chi_1^\Mid)^2x(\partial_xu_1)^2\,dx\\
&\quad\lesssim\|(x-t\partial_x^2)u_1\|_{L^2}^2+2t\int\partial_x[(\chi_1^\Mid)^2]u_1^2\,dx+t\int\partial_x^2[(\chi_1^\Mid)^2]xu_1^2\,dx\\
&\quad\lesssim \|u_1\|_{\XEquiv}^2.
\end{align*}
Using the localization we obtain
\begin{align*}
\|(|x|+t)\chi_1^\Mid u_1\|_{L^2}^2&\lesssim \|\chi_1^\Mid x u_1\|_{L^2}^2+\|\chi_1^\Mid t\partial_x^2u_1\|_{L^2}^2+\|u_1\|_{\XEquiv}^2+2t\int (\chi_1^\Mid)^2x(\partial_xu_1)^2\,dx\\
&\lesssim\|u_1\|_{\XEquiv}^2,
\end{align*}
which completes the proof of \eqref{est:Elliptic}.

\textit{C. Hyperbolic region.} We note that \(u_{N,-}^\hyp=\overline{u_{N,+}^\hyp}\) so it suffices to consider positive frequencies. By scaling it again suffices to consider the case \(N=1\). We define \(f_{1,+}=L_+u_{1,+}^\hyp\) and may argue as in the inequality \eqref{est:LocalisationHelper} to obtain
\est{KeyWjp}{
\|(1-P_{\frac14\leq\cdot\leq 4}P_+)\partial_x^\alpha(|x|^\beta f_{1,+})\|_{L^2}\lesssim_k t^{-k}\|u_1\|_{\XEquiv}.
}
We observe that
\[
\||x|^{\frac12}f_{1,+}\|_{L^2}^2+t\|\partial_xf_{1,+}\|_{L^2}^2=\|L_-f_{1,+}\|_{L^2}^2+4t\Im\int(|x|^{\frac14}f_{1,+})\partial_x\overline{(|x|^{\frac14} f_{1,+})}\,dx.
\]
Using \eqref{est:KeyWjp} we may estimate
\begin{align*}
t^{\frac12}\|f_{1,+}\|_{L^2} &\lesssim t^{\frac12}\|\partial_xf_{1,+}\|_{L^2}+\|u_1\|_{\XEquiv},\\
\|L_-f_{1,+}\|_{L^2} &\lesssim \|Lu_1\|_{L^2}+\|u_1\|_{\XEquiv},\\
4t^{\frac12}\Im\int(|x|^{\frac14}f_{1,+})\partial_x\overline {(|x|^{\frac14}f_{1,+})}\,dx &\lesssim \|u_1\|_{\XEquiv},
\end{align*}
where the last estimate uses that \(|x|^{\frac14}f_{1,+}\) is localized to positive frequencies up to rapidly decaying tails.
This gives us
\[
\|(|x|^{\frac12}+t^{\frac12})f_{1,+}\|_{L^2}^2\lesssim \||x|^{\frac12}f_{1,+}\|_{L^2}^2+t\|\partial_xf_{1,+}\|_{L^2}^2+\|u_1\|_{\XEquiv}^2\lesssim \|u_1\|_{\XEquiv}^2.
\]
\epf

\medskip
\begin{proof}[Proof of Proposition \ref{propn:InitialPointwiseBounds}]

We now turn to the proof of \eqref{est:EllipticBounds}. The \(L^2\) bound follows from the energy estimate \eqref{est:AltNorm} and the elliptic bounds \eqref{est:Low} and \eqref{est:Elliptic}. For the second part, we use Bernstein's inequality \eqref{est:Bernie} to obtain
\[
\|t^{\frac13}u^\el\|_{L^\infty(\SSR_0)}\lesssim t^{-\frac16}\|{t^{\frac13}}u_{\leq t^{-\frac13}}\|_{L^2}+\sum\limits_{N> t^{-\frac13}}t^{-\frac23}N^{-\frac32}\|t N^2 u_N^\el\|_{L^2}.
\]
Using the elliptic bound \eqref{est:Elliptic} and the Cauchy-Schwarz inequality to sum in \(N>t^{-\frac13}\) we obtain
\[
\|t^{\frac13}u^\el\|_{L^\infty(\SSR_0)}\lesssim t^{-\frac16}\|u\|_{\XEquiv}\lesssim\epsilon.
\]

For \(M>t^{-\frac13}\) we take smooth \(\chi_M\) localizing to the set \(\{|x|\sim tM^2\}\) as defined in Lemma \ref{lem:EllipticEstimates}. From \eqref{est:LocalisationHelper} \(\chi_M u_N^\el\) is localized at frequency \(\lesssim N\) for \(N\leq M\) up to rapidly decaying tails of size \(O((tM^2N)^{-k})\). As a consequence we may combine Bernstein's inequality \eqref{est:Bernie} and the inequality \eqref{est:LocalisationHelper} to obtain
\[
\|\chi_Mu^\el\|_{L^\infty}\lesssim t^{-\frac16}\|\chi_M u_{\leq t^{-\frac13}}\|_{L^2}+\sum\limits_{t^{-\frac13}<N\leq M}N^{\frac12}\|\chi_M u_N^\el\|_{L^2}+\sum\limits_{N>M}N^{\frac12}\|u_N^\el\|_{L^2}+t^{-1}M^{-\frac32}\|u\|_{\XEquiv}.
\]
We then have
\begin{align*}
\|t^{\frac13}\langle t^{-\frac13}x\rangle^{\frac34}u^\el\|_{L^\infty(|x|\sim tM^2)}&\lesssim t^{\frac56}M^{\frac32}\|\chi_Mu^\el\|_{L^\infty}\\
&\lesssim t^{-\frac13}M^{-\frac12}\|x u_{\leq t^{-\frac13}}\|_{L^2}+\sum\limits_{t^{-\frac13}<N\leq M}t^{-\frac16}N^{\frac12}M^{-\frac12}\|xu_N^\el\|_{L^2}\\
&\quad+\sum\limits_{N>M}t^{-\frac16}M^{\frac32}N^{-\frac32}\|tN^2u_N^\el\|_{L^2}+t^{-\frac16}\|u\|_{\XEquiv}.
\end{align*}
Summing in \(N\) using the Cauchy-Schwarz inequality and the energy estimate \eqref{est:AltNorm} we obtain the second part of \eqref{est:EllipticBounds} by taking the supremum over \(M>t^{-\frac13}\).

For the third part of \eqref{est:EllipticBounds} we estimate similarly for \(M>t^{-\frac13}\) to obtain
\begin{align*}
\|t^{\frac23}\langle t^{-\frac13}x\rangle^\frac14u^\el_x\|_{L^\infty(|x|\sim tM^2)}&\lesssim t^{-\frac12}M^{-\frac32}\|x u_{\leq t^{-\frac13}}\|_{L^2}+\sum\limits_{t^{-\frac13}<N\leq M}t^{-\frac16}N^{\frac32}M^{-\frac32}\|xu_N^\el\|_{L^2}\\
&\quad+\sum\limits_{N>M}t^{-\frac16}M^{\frac12}N^{-\frac12}\|tN^2u_N^\el\|_{L^2}+t^{-\frac16}\|u\|_{\XEquiv}.
\end{align*}

For \eqref{est:HyperbolicBounds} we apply the Sobolev estimate \eqref{est:SobolevEst} to \(e^{-i\phi}u_{N,+}^\hyp\) to get
\begin{align*}
\|t^{\frac13}u_{N,+}^\hyp\|_{L^\infty}&\lesssim t^{\frac{1}{12}}\|u_{N,+}^\hyp\|_{L^2}^{\frac12}\|L_+u_{N,+}^\hyp\|_{L^2}^{\frac12}\\
&\lesssim t^{-\frac16}\|t^{\frac12}NL_+u_{N,+}^\hyp\|_{L^2}+t^{-\frac16}\|N^{-1}u_N\|_{L^2}
\end{align*}
Summing over \(N\geq t^{-\frac13}\) using that the \(u_{N,\pm}^\hyp\) have almost disjoint supports and \eqref{est:Hyperbolic}, we obtain the first part of \eqref{est:HyperbolicBounds}.

For the second part we may use the localization to estimate
\[
\|t^{\frac23}\langle t^{-\frac13}x\rangle^{-\frac12}\partial_xu_{N,+}^\hyp\|_{L^\infty}\lesssim \|t^{\frac13}u_{N,+}^\hyp\|_{L^\infty}+t^{-\frac16}\|u_N\|_{\XEquiv},
\]
and then apply the first part of \eqref{est:HyperbolicBounds} to the first term.

\epf
}

\medskip

\subsection{Construction of wave packets} Let \(\chi\in C^\infty_0(\R)\) be a real-valued function, supported on a neighbourhood of the origin of size \(\sim1\) and localized in frequency near \(0\) at scale \(\sim1\) satisfying \(\int\chi=1\). We then define our wave packet
\defn{Psiv}{
\Psi_v(t,x)=\chi(\lambda(x+tv))e^{i\phi},
}
where \(\phi\), \(\lambda\) are defined as in \eqref{defn:phi}, \eqref{defn:lambda}. We define \(\VOR_\rho\) such that if \(v\in\VOR_\rho\), then \(\Psi_v\) is supported on \(\OR_\rho\). The following lemma shows that \(\Psi_v\) is also a good approximation to a free solution in Fourier space.

\medskip
\lemma{FourierTransform}{
For \(t\geq1\) and \(v\in\VOR_0\)
\eq{FSide}{
\hat\Psi_v(t,\xi)=\pi^{\frac{1}{2}}\lambda^{-1}\chi_1(\lambda^{-1}(\xi-\xiv))e^{\frac{1}{3}it\xi^3},
}
where \(\xiv=\sqrt v\), and \(\chi_1\in\Schwartz(\R)\) is localized at scale \(1\) in space and frequency satisfying
\est{IntegralBound}{
\int\chi_1(\xi)=1+O\left((t^{\frac{2}{3}}v)^{-\frac{3}{4}}\right).
}
}
\bpf
We consider the Taylor approximation of \(\phi\) at \(x=-tv\),
\[
\phi(t,x)=\tfrac{1}{3}t\xiv^3+x\xiv+\frac{\pi}{4}-\frac{1}{4}(\lambda(x+tv))^2+R(\lambda(x+tv),t^{\frac{2}{3}}v),
\]
where
\[
R(x,y)=-\int_0^1\frac{y^{-\frac{3}{4}}x^3(1-h)^2}{8|y^{-\frac{3}{4}}xh-1|^{\frac{3}{2}}}\,dh
\]
is well defined for \(x\in\supp\Psi_v\) whenever \(v\in\VOR_0\). We may then define
\[
\chi_1(\xi)=\pi^{-1}e^{-\frac{1}{3}it\lambda^3\xi^3}\int e^{-2i\xi\eta}e^{i\eta^2}\hat\chi_2(\eta)\,d\eta,
\]
where
\(
\chi_2(x)=\chi(x)e^{iR(x,t^{\frac{2}{3}}v)}.
\)

As \(e^{\frac{1}{3}it\lambda^3\xi^3}=1+O((t^{\frac{2}{3}}v)^{-\frac{3}{4}}\xi^3)\) and \(\chi_2\in\Schwartz\) we have
\[
\int\chi_1=\hat\chi_2(0)+O((t^{\frac{2}{3}}v)^{-\frac{3}{4}}),
\]
and similarly, as \(e^{iR(x,y)}=1+O(y^{-\frac{3}{4}}x^3)\),
\[
\hat\chi_2(0)=1+O((t^{\frac{2}{3}}v)^{-\frac{3}{4}}).
\]
\epf
\medskip

{
We will frequently make use of the fact that up to error terms we may replace \(u\) in the definition of \(\gamma(t,v)\) by the hyperbolic part of \(u\) lying on the ray \(\Gamma_v\). To prove this we first take \(v\in\VOR_0\) and define \(\zeta_v\in C^\infty\) such that
\[
\zeta_v(D) = \sum\limits_{N\sim \xiv}P_NP_+.
\]
As \(\Psi_v\) is localized in Fourier space at frequency \(\xiv\), from \eqref{eq:FSide} we have
\[
\|(1 - \zeta_v(D))\Psi_v\|_{L^1_x}\lesssim_k t^{\frac13}(t^{\frac23}v)^{-k}.
\]
From the initial pointwise bounds and the elliptic bounds \eqref{est:EllipticBounds} and \eqref{est:HyperbolicBounds} we obtain
\[
\left|\gamma(t,v)-\sum\limits_{N\sim \xiv} \int u_{N,+}(t,x)\overline\Psi_v(t,x)\,dx\right|\lesssim\|(1 - \zeta_v(D))\Psi_v\|_{L^1}\|u\|_{L^\infty}\lesssim_k \epsilon(t^{\frac23}v)^{-k}.
\]
We then define
\defn{udiag}{
\udiag(t,x)=e^{-i\phi}\sum\limits_{N\sim \xiv}u_{N,+}^\hyp,
}
and from the spatial localization of the hyperbolic part of \(u_{N,+}\), we have
\est{GammaThrowAway}{
\left|\gamma(t,v)-\int \udiag(t,x)\chi(\lambda(x+tv))\,dx\right|\lesssim_k \epsilon(t^{\frac23}v)^{-k}.
}
}

\subsection{Energy estimates for \(\gamma\)}

We may consider \(\gamma\) to be a function of \(\xiv=\sqrt v\), and define \(\FOR_\rho\) so that \(\xiv\in\FOR_\rho\) if and only if \(v\in\VOR_\rho\). We then have the following energy estimates for \(\gamma\).

\medskip
\lemma{GammaNRG}{
For \(t\in[1,T]\) we have the energy estimates
\begin{gather}
\|\gamma\|_{H^{0,1}_{\xiv}(\FOR_0)}\lesssim\epsilon\label{est:BasicGammaEstH01},\\
{\|\partial_{\xiv}\gamma-3t\xiv^{-1}\partial_t\gamma\|_{L^2(\FOR_0)}\lesssim\epsilon t^\delta\label{est:BasicGammaH1}}.
\end{gather}
}
\bpf
{
We first show that,
\est{Youngs}{
\left\|\int f(t,x)\chi(t^{-\frac12}\xiv^{-\frac12}(x+t\xiv^2))\,dx\right\|_{L^2_{\xiv}(\FOR_\rho)}\lesssim \|f\|_{L^2(\OR_{\rho})}.
}
Making an affine change of variables, we have
\[
\int f(t,x)\chi(t^{-\frac12}\xiv^{-\frac12}(x+t\xiv^2))\,dx=\int t^{\frac12}\xiv^{\frac12}f(t,t^{\frac12}\xi_v^{\frac12}x-t\xiv^2)\chi(x)\,dx.
\]
We then define a nonlinear change of variables by
\[
\xiv\mapsto q=t^{\frac12}\xi_v^{\frac12}x-t\xiv^2.
\]
We calculate
\[
t^{-\frac13}q=-(t^{\frac13}\xiv)^2\left(1-(t^{\frac13}\xi_v)^{-\frac32}x\right),\qquad\frac{dq}{d\xiv}=-2t\xiv\left(1-\frac14(t^{\frac13}\xi_v)^{-\frac32}x\right).
\]
If \(\xiv\in\FOR_\rho\), then \(t^{\frac13}\xiv\gtrsim t^\rho\geq1\). Provided \(\chi\) is supported in a sufficiently small neighbourhood of the origin we have
\[
-t^{-\frac13}q\gtrsim t^{2\rho},\qquad \left|\frac{dq}{d\xiv}\right|\gtrsim t\xiv,
\]
which gives us the estimate \eqref{est:Youngs}.

As a consequence of \eqref{est:Youngs},} we have the estimate
\[
\|\gamma\|_{L^2_{\xiv}(\FOR_0)}\lesssim\|u\|_{L^2}.
\]
We calculate
\[
\xiv\Psi_v=-i\partial_x\Psi_v+\lambda\tilde\Psi_v,
\]
where
\[
\tilde\Psi_v(t,x)=\left(\lambda^{-1}(\xiv-t^{-\frac{1}{2}}|x|^{\frac{1}{2}})\chi(\lambda(x+tv))+i\chi'(\lambda(x+tv))\right)e^{i\phi}
\]
has similar localization to \(\Psi_v\). Integrating by parts in the first term {and using \eqref{est:Youngs}}, we obtain
\[
\|\xiv\gamma\|_{L^2_{\xiv}(\FOR_0)}\lesssim\|u\|_{H^1}.
\]

{
We now turn to the estimate \eqref{est:BasicGammaH1}. We observe that
\(
(3t\partial_t+x\partial_x)\Psi_v = \xiv\partial_{\xiv}\Psi_v
\)
so integrating by parts we have
\[
\partial_{\xiv}\gamma - 3t\xiv^{-1}\partial_t\gamma = \int \Lambda u\,\xiv^{-1}\partial_x\overline\Psi_v\,dx.
\]
We calculate that
\[
\xiv^{-1}\partial_x\Psi_v(t,x) = \left(\xiv^{-1}\lambda\chi'(\lambda(x+tv))+it^{-\frac12}|x|^{\frac12}\xiv^{-1}\chi(\lambda(x+tv))\right)e^{i\phi}
\]
has similar localization to \(\Psi_v\). From the estimate \eqref{est:Youngs}, we then obtain
\[
\|\partial_{\xiv}\gamma - 3t\xiv^{-1}\partial_t\gamma\|_{L^2(\FOR_0)}\lesssim \|\Lambda u\|_{L^2}.
\]
}

\epf

\medskip

\subsection{Reduction of pointwise estimates to wave packets}

The following lemma allows us to reduce closing the bootstrap estimate \eqref{est:Bootstrap} to proving
\est{MainGamma}{
\|\gamma\|_{L^\infty_v(\VOR_0)}\lesssim\epsilon
}
with a constant independent of \(M_0,T\).

\medskip
\proposition{DispersiveEstimates}{
For \(t\in[1,T]\) we have the following estimates.

I) Physical-space estimates.
\begin{align}
\left\|t^{\frac{1}{3}}(t^{-\frac{1}{3}}|x|)^{\frac{3}{8}}\left(P_+u(t,x)-t^{-\frac{1}{3}}(t^{-\frac{1}{3}}|x|)^{-\frac{1}{4}}e^{i\phi}\gamma(t,t^{-1}|x|)\right)\right\|_{L^\infty(\OR_0)} &\lesssim\epsilon\label{est:OscLInf},\\
\left\|t^{\frac{2}{3}}(t^{-\frac{1}{3}}|x|)^{-\frac{1}{8}}\left(P_+u_x(t,x)-it^{-\frac{2}{3}}(t^{-\frac{1}{3}}|x|)^{\frac{1}{4}}e^{i\phi}\gamma(t,t^{-1}|x|)\right)\right\|_{L^\infty(\OR_0)} &\lesssim\epsilon\label{est:OscLInfD},\\
\left\|t^{\frac{1}{6}}(t^{-\frac{1}{3}}|x|)^{\frac{1}{4}}\left(P_+u(t,x)-t^{-\frac{1}{3}}(t^{-\frac{1}{3}}|x|)^{-\frac{1}{4}}e^{i\phi}\gamma(t,t^{-1}|x|)\right)\right\|_{L^2(\OR_0)} &\lesssim\epsilon\label{est:OscL2}.
\end{align}

II) Fourier-space estimates.
\begin{align}
\|(t^{\frac{1}{3}}\xi)^{\frac{1}{4}}(\hat u(t,\xi)-\pi^{-\frac{1}{2}}e^{\frac{1}{3}it\xi^3}\gamma(t,\xi^2))\|_{L^\infty_\xi(\FOR_0)} &\lesssim\epsilon\label{est:FouLInf},\\
\|t^{\frac{1}{6}}(t^{\frac{1}{3}}\xi)^{\frac{1}{2}}(\hat u(t,\xi)-\pi^{-\frac{1}{2}}e^{\frac{1}{3}it\xi^3}\gamma(t,\xi^2))\|_{L^2_\xi(\FOR_0)} &\lesssim\epsilon\label{est:FouL2}.
\end{align}
}

\bpf~


\textit{I) Physical-space estimates.} For \eqref{est:OscL2}, using the elliptic estimate \eqref{est:EllipticBounds} and the estimate \eqref{est:GammaThrowAway}, it suffices to show that
\[
\left\|\lambda^{-2}\udiag(t,-tv)-\lambda^{-1}\int \udiag(t,x)\chi(\lambda(x+tv))\,dx\right\|_{L^2_{\xiv}(\FOR_0)}\lesssim\epsilon t^{\frac{1}{6}}.
\]

As \(\int\chi=1\) we have
\begin{align*}
&\lambda^{-1}\udiag(t,-tv)-\int \udiag(t,x)\chi(\lambda(x+tv))\,dx\\
&\quad=\int \left(\udiag(t,-tv)-\udiag(t,x)\right)\chi(\lambda(x+tv))\,dx\\
&\quad=-\int \int_0^1(\partial_x\udiag)(t,x-(x+tv)h)(x+tv)\chi(\lambda(x+tv))\,dhdx.
\end{align*}
We observe that
\[
\partial_x \udiag = e^{-i\phi}\sum\limits_{N\sim \xiv} L_+u_{N,+}^\hyp.
\]
The estimate then follows from the hyperbolic bound \eqref{est:HyperbolicBounds} and the estimate \eqref{est:Youngs}.

{
For \eqref{est:OscLInf} we estimate similarly, using the elliptic estimate \eqref{est:EllipticBounds} and the estimate \eqref{est:GammaThrowAway} to reduce the bound to proving that
\[
\left\|\lambda^{-\frac12}\int \left(\udiag(t,-tv)-\udiag(t,x)\right)\chi(\lambda(x+tv))\,dx\right\|_{L^\infty_{\xiv}(\FOR_0)}\lesssim\epsilon t^{\frac16}.
\]
To show this we use \eqref{est:HyperbolicBounds} to obtain
\[
\lambda^{-\frac{1}{2}}|(\udiag(t,-tv)-\udiag(t,x))|\lesssim (t^{\frac{1}{3}}\xiv)^{\frac{1}{4}}\|\partial_x\udiag\|_{L^2}|x+tv|^{\frac{1}{2}}\lesssim \epsilon t^{\frac16}\lambda^{\frac{3}{2}}|x+tv|^{\frac{1}{2}}.
\]

Using \eqref{est:GammaThrowAway} and the frequency localization, the estimate \eqref{est:OscLInfD} follows from \eqref{est:OscLInf}.
}

\textit{II) Fourier-space estimates.} We use \eqref{eq:FSide} to write
\begin{align*}
e^{-\frac{1}{3}it\xiv^3}\hat u(t,\xiv)-\pi^{-\frac{1}{2}}\gamma(t,v)&=\pi^{-\frac{1}{2}}\int \left(e^{-\frac{1}{3}it\xiv^3}\hat u(t,\xiv)-e^{-\frac{1}{3}it\xi^3}\hat u(t,\xi)\right)\lambda^{-1}\chi_1(\lambda^{-1}(\xi-\xiv))\,d\xi\\
&\quad+O\left((t^{\frac{1}{3}}\xiv)^{-\frac{3}{2}}e^{-\frac{1}{3}it\xiv^3}\hat u(t,\xiv)\right).
\end{align*}
For the difference we have
\begin{align*}
e^{-\frac{1}{3}it\xiv^3}\hat u(t,\xiv)-e^{-\frac{1}{3}it\xi^3}\hat u(t,\xi)&=-i(\xiv-\xi)\int_0^1e^{-\frac{1}{3}it\eta^3}\widehat{(\J u)}(t,h(\xiv-\xi)+\xi)\,dh.
\end{align*}
{For the error terms we have
\begin{align*}
\|t^{\frac{1}{6}}(t^{\frac{1}{3}}\xiv)^{-1}e^{-\frac{1}{3}it\xiv^3}\hat u(t,\xiv)\|_{L^2_{\xiv}(\FOR_0)}\lesssim t^{-\frac16}\|t^{\frac13}\langle t^{\frac13}D_x\rangle^{-1}u\|_{L^2}.
\end{align*}
The estimate \eqref{est:FouL2} then follows from the energy estimate \eqref{est:AltNorm}.
}

For \eqref{est:FouLInf} we use that
\[
|e^{-\frac{1}{3}it\xiv^3}\hat u(t,\xiv)-e^{-\frac{1}{3}it\xi^3}\hat u(t,\xi)|\lesssim\|\J u\|_{L^2}|\xiv-\xi|^{\frac{1}{2}},
\]
and estimate similarly.
\epf

\medskip

\subsection{Global existence}


In order to prove \eqref{est:MainGamma}, we consider the ODE satisfied by \(\gamma\)
\eq{GammaODE}{
\dot \gamma(t,v)=\sigma \int (u^3)_x\overline\Psi_v\,dx+\int u\overline{(\partial_t+\tfrac{1}{3}\partial_x^3)\Psi_v}\,d\xi.
}
We then have the following estimate for \(\dot\gamma\).

\medskip
{
\lemma{ODEEsts}{
For \(t\in[1,T]\) and \(\epsilon>0\) sufficiently small, we have the estimates
\begin{align}
\|t(t^{\frac23}v)^{\frac18}(\dot\gamma - 3i\sigma t^{-1}|\gamma|^2\gamma)\|_{L^\infty_v(\VOR_0)} &\lesssim \epsilon,\label{est:GammaDotLInf}\\
\|t^{\frac76}(t^{\frac13}\xiv)^{\frac12}(\dot\gamma-3i\sigma t^{-1}|\gamma|^2\gamma)\|_{L^2_{\xiv}(\FOR_0)} &\lesssim\epsilon\label{est:GammaDotL2}.
\end{align}
}
\bpf
We use \(\err\) to denote error terms that satisfy the estimates
\[
\|t(t^{\frac23}v)^{\frac18}\err\|_{L^\infty(\VOR_0)}\lesssim\epsilon,\qquad \|t^{\frac76}(t^{\frac13}\xiv)^{\frac12}\err\|_{L^2(\FOR_0)}\lesssim\epsilon.
\]

We first integrate by parts to obtain
\[
\int (u^3)_x\overline\Psi_v\,dx = 3i\int t^{-\frac12}|x|^{\frac12} u^3\bar\Psi_v\,dx - 3\int u^3\lambda e^{-i\phi}\chi'(\lambda(x+tv))\,dx.
\]
Using the bootstrap assumption \eqref{est:Bootstrap} and elliptic estimates \eqref{est:EllipticBounds}, we then have
\[
3i\int t^{-\frac12}|x|^{\frac12} u^3\bar\Psi_v\,dx - 3\int u^3\lambda e^{-i\phi}\chi'(\lambda(x+tv))\,dx = 3i\xiv\int (u^\hyp)^3\overline\Psi_v\,dx+\err.
\]

As \(u^\hyp_{N,\pm}\) is localized at frequency \(\sim \pm N\) up to rapidly decaying tails and \(\Psi_v\) at frequency \(\sim +N\), we may estimate as in \eqref{est:GammaThrowAway} to obtain
\[
3i\xiv\int (u^\hyp)^3\overline\Psi_v\,dx = 3i\xiv\int |w_{v,+}|^2w_{v,+} \chi\,dx+\err.
\]
Estimating as in \eqref{est:OscLInf} we have
\begin{align*}
3i\xiv\int |w_{v,+}|^2w_{v,+} \chi\,dx &= 3i\xiv|w_{v,+}(t,-tv)|^2\int w_{v,+}\chi\,dx+\err\\
&= 3it^{-1}|\gamma(t,v)|^2\int w_{v,+} \chi\,dx+\err\\
&=3i t^{-1}|\gamma|^2\gamma+\err.
\end{align*}

For the linear terms we calculate,
\[
(\partial_t+\tfrac13\partial_x^3)\Psi_v=t^{-1}\lambda^{-1}e^{i\phi}\partial_x\tilde\chi-\frac14it^{-\frac12}|x|^{-\frac32}\Psi_v,
\]
where
\[
\tilde \chi = \frac12\lambda(x+tv)\chi+i\lambda^2t^{\frac12}|x|^{\frac12}\chi'+\frac13t\lambda^3\chi''
\]
has the same localization as \(\chi\). For the first of these we proceed as in \eqref{est:GammaThrowAway} to replace \(u\) by a frequency localized hyperbolic piece, integrate by parts and apply the hyperbolic bound \eqref{est:Hyperbolic} to obtain
\[
t^{-1}\lambda^{-1}\int u\overline{e^{i\phi}\partial_x\tilde\chi}\,dx = -t^{-1}\lambda^{-1}\int \partial_xw_{v,+} \overline{\tilde\chi}\,dx+\err = \err.
\]
For the second term, we may simply use the localization and the hyperbolic bound \eqref{est:HyperbolicBounds} to get
\[
\frac14i\int u t^{-\frac12}|x|^{-\frac32}\Psi_v\,dx = \err.
\]
\epf
}
\medskip


We now use Lemma \ref{lem:ODEEsts} to solve the ODE \eqref{eq:GammaODE} for \(t\in [\max\{1,Cv^{-\frac32}\},T]\), where \(C>0\) is chosen such that \(v\in\VOR_0\) for \(t\geq\max\{1,Cv^{-\frac32}\}\). For velocities \(v\geq C^{\frac23}\), the ray \(\Gamma_v\) lies outside the self-similar region for all \(t\geq1\), so from \eqref{est:LWPBnd} and \eqref{eq:FSide} we may take initial data
\est{InitHi}{
|\gamma(1,v)|\lesssim\|\hat u(1)\|_{L^\infty}\lesssim \|u(1)\|_{L^2}^\frac12\|Lu(1)\|_{L^2}^\frac12\lesssim\epsilon.
}
{
For velocities \(0<v<C^{\frac23}\), the ray \(\Gamma_v\) lies inside the self-similar region up to time \(t_0=Cv^{-\frac{3}{2}}\), so using \eqref{est:GammaThrowAway}, Bernstein's inequality \eqref{est:Bernie} and the energy estimate \eqref{est:AltNorm}, we have initial data
\est{InitLo}{
|\gamma(t_0,v)|\lesssim t_0^{\frac16}\sum\limits_{N\sim t_0^{-\frac13}}\|u_N(t_0)\|_{L^2}+\epsilon\lesssim\epsilon.
}
}

{
From \eqref{est:GammaDotLInf}, for \(v\in\VOR_0\), we have the estimate
\[
\dot\gamma=3i\sigma t^{-1}|\gamma|^2\gamma+O\left({\epsilon}t^{-1}(t^{\frac23}v)^{-\frac18}\right).
\]
As \(3i\sigma t^{-1}|\gamma|^2\) is imaginary we can then solve \eqref{eq:GammaODE} to find a solution satisfying \eqref{est:MainGamma}. This completes the proof of global existence.
}

\subsection{Asymptotic behavior}

From Lemma \ref{lem:ODEEsts} there exists a unique function \(W\) defined on \((0,\infty)\) such that for \(t\geq1\)
\begin{align}
\|(t^{\frac{2}{3}}v)^{\frac{1}{8}}(\gamma(t,v)-(2\pi)^{-\frac{1}{2}}W(\xiv)e^{\frac{3i\sigma}{4\pi}|W(\xiv)|^2{\log(t\xiv^3)}})\|_{L^\infty_v(\VOR_0)} &\lesssim\epsilon,\\
\|t^{\frac{1}{6}}(t^{\frac{1}{3}}\xiv)^{\frac{1}{2}}(\gamma(t,v)-(2\pi)^{-\frac{1}{2}}W(\xiv)e^{\frac{3i\sigma}{4\pi}|W(\xiv)|^2{\log(t\xiv^3)}})\|_{L^2_{\xiv}(\FOR_0)} &\lesssim\epsilon\label{est:L2Helper}.
\end{align}
{
We extend \(W\) to \(\R\) by defining
\[
W(-\xiv)=\overline W(\xiv),\qquad W(0)=\int u_0.
\]
}
As every \(\xiv>0\) lies in \(\FOR_0\) for sufficiently large \(t>0\), using the estimates \eqref{est:BasicGammaEstH01} and \eqref{est:MainGamma} we obtain
\est{EasyWBnds}{
\|W\|_{L^\infty_{\xiv}\cap H^{0,1}_{\xiv}}\lesssim\epsilon.
}
{
To prove the additional regularity for \(W\) we define the region \(\FOR_*=\FOR_{\sfrac12}\backslash\FOR_{\sfrac16} = \{t^{-\frac16}\lesssim\xiv\lesssim t^{\frac16}\}\) and the phase \(\Phi = 3\sigma|\gamma(t,v)|^2\log(t\xiv^3)\). From the estimates \eqref{est:L2Helper} and \eqref{est:MainGamma} we have,
\[
\left\|e^{-i\Phi}\gamma(t,v)-(2\pi)^{-\frac12}W(\xiv)\right\|_{L^2_{\xiv}(\FOR_*)}\lesssim \epsilon t^{-\frac14}(1+\log t).
\]
We calculate
\begin{align*}
e^{i\Phi}\partial_{\xiv}\left(e^{-i\Phi}\gamma(t,v)\right)&=\partial_{\xiv}\gamma - 9i\sigma\xiv^{-1} |\gamma|^2\gamma - 3i\sigma\partial_{\xiv}(|\gamma|^2)\gamma \log(t\xiv^3)\\
&=(\partial_{\xiv}\gamma-3t\xiv^{-1}\dot\gamma)-6i\sigma\Re((\partial_{\xiv}\gamma-3t\xiv^{-1}\dot\gamma)\overline\gamma)\gamma\log(t\xiv^3)\\
&\quad+O\left(t\xiv^{-1}\left|\dot\gamma-3i\sigma|\gamma|^2\gamma\right|\left(1+|\gamma|^2\log(t\xiv^3)\right)\right).
\end{align*}
From the energy estimate \eqref{est:BasicGammaH1} and the estimate \eqref{est:GammaDotL2}, we have
\[
\left\|\partial_{\xiv}\left(e^{-i\Phi}\gamma(t,v)\right)\right\|_{L^2(\FOR_*)}\lesssim \epsilon t^\delta(1+\log t).
\]
}
By interpolation, for sufficiently large \(C>0\), we then have
\est{LargepH1Bnd}{
\|W\|_{H^{1-C\epsilon^2}_{\xiv}}\lesssim \epsilon.
}

To derive the asymptotic behavior in the self-similar region, we use the self-similar change of variables \eqref{defn:BigU}.
Let \(\rho>0\) and \(C\gg1\). From the equation for \(U\) \eqref{eq:mKdVForBigU}, the energy estimate \eqref{est:JNRG}, Bernstein's inequality \eqref{est:Bernie} and the elliptic estimate \eqref{est:Elliptic}, we have
\begin{align*}
\|\partial_tP_{\leq Ct^{\rho}}U\|_{L^\infty(\SSR_\rho)}&\lesssim t^{\frac\rho2}\|P_{\leq Ct^{\rho}}\partial_tU\|_{L^2}+t^{-1}\|P_{\sim Ct^\rho}U\|_{L^\infty(\SSR_\rho)}\\
&\lesssim\epsilon t^{\frac32\rho+\delta-\frac76}+\epsilon t^{\frac\rho2-\frac56}\sum\limits_{N\sim Ct^{\rho-\frac13}}\|u_N^\el\|_{L^2}\\
&\lesssim \epsilon t^{-\min\{\frac16-\delta-\frac32\rho,\frac32\rho\}-1}.
\end{align*}
{Further, from the elliptic estimate \eqref{est:Elliptic} we also have
\[
\|P_{>Ct^{\rho}}U\|_{L^\infty(\SSR_\rho)}\lesssim \sum\limits_{N>Ct^{\rho-\frac13}}t^{\frac16}N^{\frac12}\|u_N^\el\|_{L^2}\lesssim \epsilon t^{-\frac32\rho}.
\]
Choosing \(0<\rho<\frac23(\frac16-\delta)\) there exists \(Q\in L^\infty\) such that
\[
\|Q\|_{L^\infty}\lesssim\epsilon,\qquad\|U-Q\|_{L^\infty(\SSR_\rho)}\lesssim \epsilon t^{-\min\{\frac16-\delta-\frac32\rho,\frac32\rho\}}. 
\]

We recall that
\[
\|yU-U_{yy}+3\sigma U^3\|_{L^2_y}\lesssim\epsilon t^{\delta-\frac16},
\]
and taking the limit as \(t\rightarrow\infty\), we have that \(Q\) is a solution to the Painlev\'e II equation \eqref{eq:PLII}.
}

%
%
\section{Asymptotic Completeness}\label{sect:AsympComplete}
In this section we prove Theorem \ref{thrm:AC}. We note that from the local theory it suffices to prove existence of a solution \(u(t)\) on \([1,\infty)\) satisfying
\est{v1Bnd}{
\|u(1)\|_{\X}\lesssim\epsilon.
}
{For \(C,\epsilon\) as in \eqref{est:WAC}, we define \(\delta=C\epsilon^2\).}

\subsection{Regularization of \(W\)} Instead of working with \(\uasymp\), we work with an approximation \(\uapp\) given by regularizing \(W\) at the scale corresponding to the wave packets.

{We dyadically decompose
\[
W(z)=\sum\limits_{N\in 2^\Z}W_N(z),\qquad W_N=P_NW.
\]
Let \(\chi\in C^\infty\) be smooth on scale \(\sim1\), such that \(\chi(z)\equiv1\) for \(|z|\geq1\) and \(\chi(z)\equiv0\) for \(|z|\leq\frac12\). For each \(N>1\) we define the function
\[
\chi_N(t,z)=\chi(N^{-2}t^{\frac23}\langle t^{\frac{1}{3}}z\rangle),
\]
We observe that \(\chi_N\equiv1\) for \(N\leq t^{\frac13}\) and for \(N>t^{\frac13}\) they are localized on the set \(A_N=\{t^{\frac{2}{3}}\langle t^{\frac{1}{3}}z\rangle\gtrsim N^2\}\) and at frequencies \(\lesssim tN^{-2}< N\) up to rapidly decaying tails.} We then define
\defn{WReg}{
\W(t,z)=\sum\limits_{N\leq t}\chi_N(t,z)W_N(z).
}
By construction, the map
\(
x\mapsto \W(t,t^{-\frac{1}{2}}|x|^{\frac{1}{2}})
\)
is smooth on the scale of the wave packets on \(\R\backslash\{0\}\). However, to ensure that \(\uapp\) is a good approximation on \(\R\) we require additional smoothing at \(x=0\). {To do this we take {an even function} \(\zeta\in C^\infty\) so that \(\zeta(y)=|y|^\frac12\) for \(|y|\geq1\) and \(\zeta(y)=\frac{1}{4}\langle 16y\rangle^{\frac12}\) for \(|y|\leq \frac12\). Taking \(Q(y;W)\) as in Theorem \ref{thrm:PLII} we then define the corresponding approximate solution to be
\defn{uapp}{
\uapp(t,x)=t^{-\frac{1}{3}}Q\left(t^{-\frac{1}{3}}x;\W\left(t,t^{-\frac{1}{3}}\zeta(t^{-\frac13}x)\right)\right).
}
}

%
As a straightforward consequence of the localization of \(\W\) we have the following Lemma.

\medskip
\lemma{WEsts}{
For \(t\geq1\) we have the following estimates.

I) Estimates for \(\W=\W(t,t^{-\frac{1}{3}}\zeta(t^{-\frac13}x))\).
\est{WL2}{
\begin{gathered}
\|t^{-\frac{1}{3}}\langle t^{-\frac{1}{3}}x\rangle^{-\frac{1}{4}}\W\|_{L^2}\lesssim\epsilon,\qquad\|t^{\frac{1}{3}}\langle t^{-\frac{1}{3}}x\rangle^{\frac{1}{4}}\partial_x\W\|_{L^2}\lesssim\epsilon,\\
\begin{aligned}
\|(t^{\frac{1}{3}}\langle t^{-\frac{1}{3}}x\rangle^{\frac{1}{4}})^{k+\delta}\partial_x^k\W\|_{L^2} &\lesssim\epsilon,\quad k\geq2,\\
\|t^{\frac{1}{3}}\langle t^{-\frac{1}{3}}x\rangle^{\frac{1}{4}}\log\langle t^{-\frac13}x\rangle\partial_x\W\|_{L^2} &\lesssim \delta^{-1}\epsilon(1+\epsilon^2\log t).
\end{aligned}
\end{gathered}
}
\est{WH01}{
\|t^{-1}(t^{\frac{1}{3}}\langle t^{-\frac{1}{3}}x\rangle^{\frac{1}{4}})\W\|_{L^2} \lesssim\epsilon,\qquad
\|t^{-1}(t^{\frac{1}{3}}\langle t^{-\frac{1}{3}}x\rangle^{\frac{1}{4}})^{k+1+\delta}\partial_x^k\W\|_{L^2} \lesssim\epsilon,\quad k\geq1.
}
\est{WLInf}{
\|\W\|_{L^\infty}\lesssim\epsilon,\qquad
\|(t^{\frac{1}{3}}\langle t^{-\frac{1}{3}}x\rangle^{\frac{1}{4}})^{k+\frac{1}{2}+\delta}\partial_x^k\W\|_{L^\infty}\lesssim\epsilon,\quad k\geq1.
}

II) Estimates for \(\W_t=\W_t(t,t^{-\frac{1}{3}}\zeta(t^{-\frac13}x))\). For \(k\geq0\),
\est{WtL2}{
\|t(t^{\frac{1}{3}}\langle t^{-\frac{1}{3}}x\rangle^{\frac{1}{4}})^{k+\delta}\partial_x^k\W_t\|_{L^2} \lesssim\epsilon,\qquad \|(t^{\frac{1}{3}}\langle t^{-\frac{1}{3}}x\rangle^{\frac{1}{4}})^{k+1+\delta}\partial_x^k\W_t\|_{L^2} \lesssim\epsilon.
}

III) Estimates for \(W-\W\).
\est{WDiffL2}{
\|(t^{\frac{1}{3}}\langle t^{-\frac{1}{3}}x\rangle^{\frac{1}{4}})^{\delta}(W-\W)\|_{L^2} \lesssim\epsilon,\qquad \|(t^{\frac{1}{3}}\langle t^{-\frac{1}{3}}x\rangle^{\frac{1}{4}})^{\frac{1}{2}+\delta}(W-\W)\|_{L^\infty} \lesssim\epsilon.
}
}
\bpf

We consider the regions \(\OR_0\cup\DR_0\) and \(\SSR_0\) separately.

For {\(|y|\gg 1\)} we have \(\zeta(y)=|y|^{\frac12}\), so by changing variables we obtain
\begin{gather*}
\|t^{-\frac13}\langle t^{-\frac13}x\rangle^{-\frac14}\W\|_{L^2({\OR_0\cup\DR_0})}\lesssim\|W\|_{L^2},\\
\|t^{\frac13}\langle t^{-\frac13}x\rangle^{\frac14}\partial_x\W\|_{L^2({\OR_0\cup\DR_0})}\lesssim\|W\|_{H^1}.
\end{gather*}

Next we consider
\begin{align*}
\|t^{\frac13}\langle t^{-\frac13}x\rangle^{\frac14}\log\langle t^{-\frac13}x\rangle\partial_x\W\|_{L^2({\OR_0\cup\DR_0})} &\lesssim \|W\|_{H^1}(1 + \log t) + \|\log\langle z\rangle \partial_zW\|_{L^2}\\
&\lesssim\|W\|_{H^1}(1 + \log t) + \delta^{-1}\|\langle z\rangle^\delta \partial_zW\|_{L^2}.
\end{align*}
We may then estimate \(\|\langle z\rangle^\delta \partial_zW\|_{L^2}\lesssim \|W\|_Y\) by interpolation.

For \(k\geq2\) we calculate
\[
\partial_x^k\W=\sum\limits_{m=1}^kc_{m,k}t^{-\frac{k+m}{3}}(t^{-\frac{1}{3}}|x|)^{\frac{m}{2}-k}(\partial_z^m\W)(t,t^{-\frac{1}{2}}|x|^{\frac{1}{2}}),
\]
so after a change of variables we have
\[
\|(t^{\frac{1}{3}}\langle t^{-\frac{1}{3}}x\rangle^{\frac{1}{4}})^{k+\delta}\partial_x^k\W\|_{L^2_x({\OR_0\cup\DR_0})}\lesssim\sum\limits_{m=1}^k\|t^{\frac{1+\delta-m}{3}}(t^{\frac{1}{3}}z)^{\frac{1+\delta+2m-3k}{2}}\partial_z^m\W\|_{L^2_z(\FOR_0)}.
\]
{For \(m=k\) we have
\begin{align*}
&\|(t^{\frac{1}{3}}(t^{\frac{1}{3}}z)^{\frac12})^{1+\delta-k}\partial_z^kW_{\leq t^{\frac13}}\|_{L^2(\FOR_0)}^2+\sum\limits_{N>t^{\frac13}}\|(t^{\frac{1}{3}}(t^{\frac{1}{3}}z)^{\frac12})^{1+\delta-k}\partial_z^k(\chi_NW_N)\|_{L^2}^2\\
&\lesssim t^{\frac{2(1+\delta-k)}{3}}\|\partial_z^kW_{\leq t^{\frac13}}\|_{L^2}^2+\sum\limits_{N>t^{\frac13}}N^{2(1+\delta)}\|W_N\|_{L^2}^2\\
&\lesssim\|W\|_{H^{1+\delta}}^2.
\end{align*}
For \(1\leq m<k\) we estimate \(\W\) in \(L^\infty\) and use Bernstein's inequality to obtain,
\begin{align*}
\|t^{\frac{1+\delta-m}{3}}(t^{\frac{1}{3}}z)^{\frac{1+\delta+2m-3k}{2}}\partial_z^m\W\|_{L^2_z(\FOR_0)}^2 &\lesssim \|t^{\frac{1+\delta-m}{3}}(t^{\frac{1}{3}}z)^{\frac{1+\delta+2m-3k}{2}}\|_{L^2(|z|\gtrsim t^{-\frac13})}^2\|\partial_z^mW_{\leq t^{\frac13}}\|_{L^\infty}^2\\
&\quad +\sum\limits_{N>t^{\frac13}}\|t^{\frac{1+\delta-m}{3}}(t^{\frac{1}{3}}z)^{\frac{1+\delta+2m-3k}{2}}\|_{L^2(|z|\gtrsim t^{-1}N^2)}^2\|\partial_z^m(\chi_NW_N)\|_{L^\infty_z}^2\\
&\lesssim \|W\|_{H^{1+\delta}}^2.
\end{align*}
}
The remaining \(L^2\)-estimates \eqref{est:WH01} and \eqref{est:WtL2} in {\(\OR_0\cup\DR_0\)} are similar.

{Next we consider \(x\in\SSR_0\). Here we have
\[
\partial_x^k\W=\sum\limits_{m=1}^k c_{m,k}t^{-\frac13(k+m)} R(t^{-\frac13}x) \partial_z^m\W
\]
where \(R\) is a smooth, bounded function depending on \(\zeta\). Applying the the Cauchy-Schwarz inequality and Bernstein inequalities we obtain
\begin{gather*}
\|t^{-\frac13}\W\|_{L^2(\SSR_0)}\lesssim {t^{-\frac16}\|W_{\lesssim t^{\frac13}}\|_{L^\infty}\lesssim} \|W\|_{L^2},\qquad\|t^{\frac13}\partial_x\W\|_{L^2(\SSR_0)}\lesssim {t^{-\frac16}\|\partial_zW_{\lesssim t^{\frac13}}\|_{L^\infty}\lesssim} \|W\|_{H^1},\\
\|t^{\frac13(k+\delta)}\partial_x^k\W\|_{L^2(\SSR_0)}\lesssim {\sum\limits_{m=1}^k t^{\frac13(\frac12+\delta-m)}\|\partial_z^mW_{\lesssim t^{\frac13}}\|_{L^\infty}\lesssim}\|W\|_{H^{1+\delta}}
\end{gather*}
}
The \(L^2\)-estimates \eqref{est:WH01} and \eqref{est:WtL2} in \(\SSR_0\) are similar.

Next we turn to the \(L^\infty\) estimate \eqref{est:WLInf}. For the undifferentiated term we may simply use Sobolev embedding to obtain \(\|\W\|_{L^\infty}\lesssim \|W\|_{H^1}\). For the second part we first observe that
\begin{align*}
\|(t^{\frac13}\langle t^{-\frac13} x\rangle^{\frac14})^{k+\frac12+\delta}\partial_x^k\W\|_{L^\infty} &\lesssim \|(t^\frac13\langle t^{-\frac13} x\rangle^{\frac14})^{k+\frac12+\delta}\partial_x^k\W\|_{L^\infty(\SSR_0)}\\
&\quad + \sup\limits_{M>t^\frac13}\|(t^\frac13\langle t^{-\frac13} x\rangle^{\frac14})^{k+\frac12+\delta}\partial_x^k\W\|_{L^\infty(|x|\sim t^{-1}M^4)}.
\end{align*}
For the self-similar region \(\SSR_0\) we may use Bernstein's inequality to obtain
\[
\|(t^\frac13\langle t^{-\frac13} x\rangle^{\frac14})^{k+\frac12+\delta}\partial_x^k\W\|_{L^\infty(\SSR_0)}\lesssim \sum\limits_{m=1}^kt^{\frac16+\frac\delta3 - \frac m3}\|\partial_z^mW_{\leq t^{\frac13}}\|_{L^\infty}\lesssim \|W_{\leq t^{\frac13}}\|_{H^{1+\delta}}.
\]
In the region \(\OR_0\cup\DR_0\) we consider each set \(\{|x|\sim t^{-1}M^4\}\) separately. Using the spatial localization of the \(\chi_N\) and Bernstein's inequality we obtain
\[
\|(t^\frac13\langle t^{-\frac13} x\rangle^{\frac14})^{k+\frac12+\delta}\partial_x^k\W\|_{L^\infty(|x|\sim t^{-1}M^4)} \lesssim \sum\limits_{k=1}^mt^{k-m}M^{2m-3k+\frac12+\delta}\|\partial_z^mW_{\lesssim M}\|_{L^\infty}\lesssim \|W\|_{H^{1+\delta}}.
\]

For the estimate for the difference \eqref{est:WDiffL2} we write
\[
W-\W=\sum\limits_{{t^{\frac13}<N<t}}(1-\chi_N)W_N + W_{>t}.
\]
For the first term we may simply estimate as above using that \(1-\chi_N\) is localized on the complement of \(A_N\). For the second term we have
\begin{align*}
\|(t^{\frac13}\langle t^{-\frac13}x\rangle^{\frac14})^\delta W_{>t}\|_{L^2} &\lesssim t^{1+\delta}\|W_{>t}\|_{L^2(|z|\leq t)} + t^\delta\|zW_{>t}\|_{L^2(|z|>t)}\\
&\lesssim \|W\|_{H^{1+\delta}} + \|\langle z\rangle \langle D\rangle^\delta W\|_{L^2}
\end{align*}

\epf

\medskip

\subsection{Estimates for \(\uapp\)} We now look to derive estimates for \(\uapp\). We first state the following lemma giving estimates for solutions to the Painlev\'e II equation \eqref{eq:PLII}, which can be proved using variation of parameters and arguments {similar to \cite{MR2309566}.}

\medskip
\lemma{PLIIEsts}{
Let \(W\in\R\) satisfy \(|W|\ll1\) and let \(Q(y;W)\) be the solution to \eqref{eq:PLII} satisfying \eqref{est:PLIIAsymp-}. We then have the estimate
\est{PLIIDecay}{
|\partial_y^k\partial_w^mQ(y;W)|\lesssim_{k,m} \begin{cases}|W|\langle y\rangle^{-\frac14+\frac k2}e^{-\frac23y_+^{\frac32}}(1+|W|^2\log \langle y\rangle)^m,&\qquad m\textrm{ even},\vspace{0.3cm}\\
\langle y\rangle^{-\frac14+\frac k2}e^{-\frac23y_+^{\frac32}}(1+|W|^2\log \langle y\rangle)^m,&\qquad m\textrm{ odd}.
\end{cases}
}

In particular, if \(|W|\leq \epsilon\) we have the estimate
\est{PLIIDecay2}{
|\partial_y^k\partial_w^mQ(y;W)|\lesssim_{k,m} \begin{cases}|W|\langle y\rangle^{-\frac14+\frac k2 +\frac\delta4}e^{-\frac23y_+^{\frac32}},&\qquad m\textrm{ even},\vspace{0.3cm}\\
\langle y\rangle^{-\frac14+\frac k2 + \frac\delta4}e^{-\frac23y_+^{\frac32}},&\qquad m\textrm{ odd}.
\end{cases}
}
}
\medskip


Using the estimates of Lemmas \ref{lem:WEsts} and \ref{lem:PLIIEsts} we can show that \(\uapp\) is a good approximation to \(\uasymp\).

\medskip
\lemma{uappGoodApprox}{
For \(t\geq1\) we have estimates for \(\uapp\)
\begin{gather}
\|t^{\frac{1}{3}}\langle t^{-\frac{1}{3}}x\rangle^{\frac{1}{4}}e^{\frac{2}{3}t^{-\frac{1}{2}}x_+^{\frac{3}{2}}}\uapp\|_{L^\infty}\lesssim\epsilon,\qquad{\|t^{\frac{2}{3}}\langle t^{-\frac{1}{3}}x\rangle^{-\frac{1}{4}}e^{\frac{2}{3}t^{-\frac{1}{2}}x_+^{\frac{3}{2}}}(\uapp)_x\|_{L^\infty}\lesssim\epsilon,}\label{est:uappBndsInf}\\
\|\uapp\|_{H^1}\lesssim\epsilon,\qquad\|L\uapp+3\sigma t\uapp^3\|_{L^2}\lesssim\epsilon(1+\epsilon^2\log t)\label{est:uappBndsNRG},
\end{gather}
and estimates for the difference \(\uapp-\uasymp\)
\begin{gather}
\|t^{\frac{1+\delta}{3}}\langle t^{-\frac{1}{3}}x\rangle^{\frac{1}{4}}(\uapp-\uasymp)\|_{L^2_x}\lesssim\epsilon\label{est:uappGoodApproxL2},\\
\|t^{\frac12+\frac\delta3}\langle t^{-\frac{1}{3}}x\rangle^{\frac{3}{8}}e^{\frac{2}{3}t^{-\frac{1}{2}}x_+^{\frac{3}{2}}}\left(\uapp-\uasymp\right)\|_{L^\infty_x}\lesssim\epsilon\label{est:uappGoodApproxLinf}.
\end{gather}

Further, if \(T\geq1\) is a dyadic integer we have the estimate
\begin{gather}
\|\uapp\|_{L^4_xL^\infty_T}\lesssim\epsilon T^{-\frac14}\label{est:uappL4LInf},
\end{gather}
where we use the notation \(L^p_T=L^p([T,2T])\).
}
\bpf
{For the first part of \eqref{est:uappBndsInf} we use the estimates \eqref{est:WLInf} for \(W\) and \eqref{est:PLIIDecay} for \(Q\) to obtain
\[
\|t^{\frac{1}{3}}\langle t^{-\frac{1}{3}}x\rangle^{\frac{1}{4}}e^{\frac{2}{3}t^{-\frac{1}{2}}x_+^{\frac{3}{2}}}\uapp\|_{L^\infty}\lesssim \|\W\|_{L^\infty}\lesssim\epsilon.
\]
For the second part we differentiate,
\[
\partial_x\uapp=t^{-\frac{2}{3}}Q_y(t^{-\frac{1}{3}}x;\W)+t^{-\frac{1}{3}}Q_w(t^{-\frac{1}{3}}x;\W)\partial_x\W,
\]
and estimate similarly,
\[
\|t^{\frac{2}{3}}\langle t^{-\frac{1}{3}}x\rangle^{-\frac{1}{4}}e^{\frac{2}{3}t^{-\frac{1}{2}}x_+^{\frac{3}{2}}}(\uapp)_x\|_{L^\infty}\lesssim \|\W\|_{L^\infty} + \|t^{\frac13}\langle t^{-\frac13}x\rangle^{-\frac12+\frac\delta4}\partial_x\W\|_{L^\infty}\lesssim\epsilon.
\]

For the first part of \eqref{est:uappBndsNRG} we estimate similarly using the estimates \eqref{est:WL2} and \eqref{est:WLInf} for \(\W\) and \eqref{est:PLIIDecay} for \(Q\) to obtain
\begin{gather*}
\|\uapp\|_{L^2}\lesssim \|t^{-\frac13}\langle t^{-\frac13}x\rangle^{-\frac14}\W\|_{L^2}\lesssim\epsilon\\
\|(\uapp)_x\|_{L^2}\lesssim t^{-\frac23}\langle t^{-\frac13} x\rangle^{\frac14}\W\|_{L^2} + \|t^{-\frac13}\langle t^{-\frac13}x\rangle^{-\frac14+\frac\delta4}\partial_x\W\|_{L^2}\lesssim \epsilon + \epsilon t^{-\frac23}.
\end{gather*}
For the second part} we use that \(Q\) satisfies \eqref{eq:PLII} to write
\[
L\uapp+3\sigma t\uapp^3=-2t^{\frac13}Q_{wy}\partial_x\W-t^{\frac23}Q_w\partial_x^2\W-t^{\frac23}Q_{ww}(\partial_x\W)^2,
\]
and may then estimate similarly. We note that the logarithmic loss arises from the first term,
\[
\|t^{\frac13}Q_{wy}\partial_x\W\|_{L^2}\lesssim \|t^{\frac13}\langle t^{-\frac13}x\rangle^{\frac14}\partial_x\W\|_{L^2} + \|\W\|_{L^\infty}^2\|t^{\frac13}\langle t^{-\frac13}x\rangle^{\frac14}\log\langle t^{-\frac13}x\rangle \partial_x\W\|_{L^2}\lesssim \epsilon(1+\epsilon^2\log t)
\]

For \eqref{est:uappGoodApproxL2} and \eqref{est:uappGoodApproxLinf} we write the difference as
\[
\uapp-\uasymp=\int_0^1t^{-\frac{1}{3}}Q_w(t^{-\frac{1}{3}}x;h\W+(1-h)W)(\W-W)\,dh,
\]
and estimate similarly using the estimates \eqref{est:WDiffL2} for the difference and \eqref{est:PLIIDecay} for \(Q\).

{
To prove \eqref{est:uappL4LInf} we take a dyadic partition of unity
\(
1=\sum_{M\in2^\Z}\varphi_M^2
\)
and taking \(l^p\) to correspond to summation in \(M\),
\begin{align*}
\|\uapp\|_{L^4_xL^\infty_T}&\lesssim\left(\sum\limits_M\|\varphi_M(t^{-\frac13}\langle t^{-\frac13}x\rangle^{\frac12})|\uapp|\|_{L^4_xL^\infty_T}^2\right)^{\frac12}\\
&\lesssim \|T^{-\frac13}\langle T^{-\frac13}x\rangle^{-\frac14}\|_{l^\infty L^4_x}\|\W\|_{l^2L^\infty_{T,x}}\\
&\lesssim T^{-\frac14}\|W\|_{H^1},
\end{align*}
where the last line follows from Sobolev embedding \eqref{est:SobolevEst} and the Cauchy-Schwarz inequality.
}

\epf

\medskip

\subsection{An equation for \(v=u-\uapp\)}


We now define the function \(f\) such that \(\uapp\) satisfies the equation
\eq{uapp}{
(\partial_t+\tfrac{1}{3}\partial_x^3)\uapp=\sigma(\uapp^3)_x+f.
}
If we define \(v=u-\uapp\) then \eqref{eq:ACPDE} becomes
\eq{ACEqnNew}{
\pde{
(\partial_t+\tfrac{1}{3}\partial_x^3)v=\mathbf N(\uapp,v)-f,
}{
\lim\limits_{t\rightarrow+\infty}v(t)=0,
}
}
where
\[
\mathbf N(\uapp,v)=\sigma\left((v+\uapp)^3-\uapp^3\right)_x.
\]

We define the norms
\begin{gather*}
\|u\|_\K=\sup\limits_{T\geq1}\left\{T^{\frac13+\frac{\delta}{3}}\|u\|_{L^\infty_TL^2_x}+T^{\frac14+\frac\delta3}\|u\|_{L^4_xL^\infty_T}+T^{\frac{\delta}{3}}\|u_x\|_{L^\infty_TL^2_x}\right\},\\
\|u\|_{\tilde\K}=\sup\limits_{T\geq1}\left\{\frac{T^{\frac{\delta}{3}}}{1+\epsilon^2\log T}\|u\|_{L^\infty_TL^2_x}\right\},
\end{gather*}
and look to solve \eqref{eq:ACEqnNew} using a contraction mapping argument in the {space
\defn{Kepsilon}{
\K_\epsilon=\{v:\|v\|_\K+\|Lv\|_{\tilde\K}\leq B\epsilon\}.
}
Instead of working with \(Lv\), we will again work with a modification
\[
\Gamma v=Lv+3\sigma t\left((v+\uapp)^3-\uapp^3\right).
\]
}

If we define \(\Phi\) as in \eqref{defn:BigPhi} then the solution to \eqref{eq:ACEqnNew} satisfies
\[
v=\Phi \mathbf N-\Phi f,\qquad\Gamma v= \Phi\tilde{\mathbf N}-\Phi\tilde f,
\]
where
\begin{gather*}
\tilde{\mathbf N}=3\sigma(v+\uapp)^2(\Gamma v)_x+3\sigma(v^2+2v\uapp)(L\uapp+3\sigma t\uapp^3)_x,\\
\tilde f=\J f+9\sigma t\uapp^2f.
\end{gather*}


\subsection{Nonlinear estimates}

For the nonlinear term we have the following estimates.

\medskip
\lemma{NonlinearEst}{
Let \(T\geq1\) be a dyadic integer and \(v_1,v_2\in \K_\epsilon\) where \(\K_\epsilon\) is defined as in \eqref{defn:Kepsilon}. Then, if \(\delta\) is defined as in \eqref{defn:delta}, for \(M_0>0\) chosen sufficiently large and \(\epsilon>0\) chosen sufficiently small we have the estimates
\begin{align}
\|v_1-v_2\|_\K+\|Lv_1-Lv_2\|_{\tilde\K} &\sim \|v_1-v_2\|_\K+\|\Gamma v_1-\Gamma v_2\|_{\tilde \K}\label{est:NonlinearHelper},\\
\|\Phi(\mathbf N(\uapp,v_1)-\mathbf N(\uapp,v_2))\|_{\K} &\ll \|v_1-v_2\|_\K\label{est:NonlinearH1}\\
\|\Phi(\tilde{\mathbf N}(\uapp,v_1)-\tilde{\mathbf N}(\uapp,v_2))\|_{\tilde\K} &\ll \|v_1-v_2\|_\K+\|Lv_1-Lv_2\|_{\tilde\K}\label{est:NonlinearH01}
\end{align}
}
\bpf
It suffices to consider \(v_1=v\), \(v_2=0\) as the general case follows by applying identical estimates. {

We first note that from the {bounds for the Airy function, we have
\[
|v|\lesssim t^{-\frac13}\langle t^{-\frac13} x\rangle^{-\frac14}\|\langle t^{-\frac13}x\rangle^{\frac14}S(-t)v\|_{L^1_x},\qquad |v_x|\lesssim t^{-\frac23}\langle t^{-\frac13}x\rangle^{\frac14}\|\langle t^{-\frac13}x\rangle^{\frac14}S(-t)v\|_{L^1_x},
\]
and hence we may estimate
\[
\|\langle t^{-\frac13}x\rangle^{\frac14}S(-t)v\|_{L^\infty_TL^1_x}\lesssim T^{-\frac\delta3}(\|v\|_Z+\|Lv\|_{\tilde Z}).
\]

Using the estimate \eqref{est:uappBndsInf}, we have
\begin{align*}
\|(v+\uapp)^3-\uapp^3\|_{L^\infty_TL^2_x}&\lesssim (\|v\|_{L^\infty_{T,x}}+\|\uapp\|_{L^\infty_{T,x}})^2\|v\|_{L^\infty_TL^2_x}\\
&\lesssim \epsilon^2T^{-1-\frac\delta3}\|v\|_{\K},
\end{align*}
Choosing \(\epsilon>0\) sufficiently small, we obtain \eqref{est:NonlinearHelper}.
}

From the estimate \eqref{est:KPV} we have
\[
\|\Phi \mathbf N\|_{\K}\lesssim\sup\limits_{T_0\geq1}\left(T_0^{\frac{1+\delta}{3}}\sum\limits_{T\geq T_0}\|\partial_x^{-1}\mathbf N\|_{L^1_xL^2_T}+T_0^{\frac\delta3}\sum\limits_{T\geq T_0}\|\mathbf N\|_{L^1_xL^2_T}\right),
\]
where we assume \(T,T_0\) are dyadic integers.

Using the estimate \eqref{est:uappL4LInf} for \(\uapp\) we may estimate
\begin{align*}
\|\partial_x^{-1}\mathbf N\|_{L^1_xL^2_T}&\lesssim (\|v\|_{L^4_xL^\infty_T}+\|\uapp\|_{L^4_xL^\infty_T})^2\|v\|_{L^2_TL^2_x}\\
&\lesssim T^{-\frac{1+\delta}{3}}(T^{-\frac\delta3}\|v\|_{\K}+\epsilon)^2\|v\|_{\K},
\end{align*}
and similarly, using the estimates \eqref{est:uappBndsNRG} and \eqref{est:uappL4LInf} for \(\uapp\),
\[
\|\mathbf N\|_{L^1_xL^2_T}\lesssim T^{-\frac{\delta}{3}}(T^{-\frac\delta3}\|v\|_{\K}+\epsilon)^2\|v\|_{\K}.
\]
Summing over dyadic \(T\geq T_0\) and using that \(\delta^{-1}\epsilon^2\lesssim M_0^{-1}\ll1\) we have \eqref{est:NonlinearH1}.

{
For \(\Phi\tilde{\mathbf N}\) we decompose
\(
\tilde{\mathbf N} = \tilde{\mathbf N}_1 - \tilde{\mathbf N}_2,
\)
where
\begin{align*}
\tilde{\mathbf N}_1 &=\partial_x\left(3\sigma(v+\uapp)^2\Gamma v+3\sigma(v^2+2v\uapp)(L\uapp+3\sigma t\uapp^3)\right),\\
\tilde{\mathbf N}_2 &= 6\sigma(v+\uapp)(v+\uapp)_x\Gamma v+3\sigma(v^2+2v\uapp)_x(L\uapp+3\sigma t\uapp^3).
\end{align*}

We may then estimate \(\tilde{\mathbf N}_1\) as before to get
\begin{align*}
\|\partial_x^{-1}\tilde{\mathbf N}_1\|_{L^1_xL^2_T} &\lesssim T^{\frac12}(\|v\|_{L^4_xL^\infty_T}+\|\uapp\|_{L^4_xL^\infty_T})^2\|\Gamma v\|_{L^\infty_TL^2_x}\\
&\quad +T^{\frac12}\|v\|_{L^4_xL^\infty_T}(\|v\|_{L^4_xL^\infty_T}+\|\uapp\|_{L^4_xL^\infty_T})\|L\uapp+3\sigma t\uapp^3\|_{L^\infty_TL^2_x}\\
&\lesssim T^{-\frac\delta3}(1+\epsilon^2\log T)(T^{-\frac\delta3}\|v\|_Z+\epsilon)^2\|\Gamma v\|_{\tilde Z}+\epsilon T^{-\frac\delta3}(1+\epsilon^2\log T)(T^{-\frac\delta3}\|v\|_Z+\epsilon)\|v\|_{Z}.
\end{align*}

For \(\tilde{\mathbf N}_2\), we use the dispersive estimates to obtain
\begin{align*}
\|\tilde{\mathbf N}_2\|_{L^1_TL^2_x} &\lesssim T\|(v+\uapp)(v+\uapp)_x\|_{L^\infty_{T,x}}\|\Gamma v\|_{L^\infty_TL^2_x}\\
&\quad +T\|(v^2+2v\uapp)_x\|_{L^\infty_{T,x}}\|L\uapp+3\sigma t\uapp^3\|_{L^\infty_TL^2_x}\\
&\lesssim T^{-\frac\delta3}(1+\epsilon^2\log T)(\|v\|_Z+\|Lv\|_{\tilde Z}+\epsilon)^2\|\Gamma v\|_{\tilde Z}+\epsilon T^{-\frac\delta3}(1+\epsilon^2\log T)\|v\|_Z(\|v\|_Z+\|Lv\|_{\tilde Z}+\epsilon)
\end{align*}

The estimate \eqref{est:NonlinearH01} then follows from applying the estimate \eqref{est:KPV} with the bounds for \(\mathbf N_1,\mathbf N_2\).

}
}
\epf
\medskip

\subsection{Inhomogeneous estimates}


In order to complete the proof of Theorem \ref{thrm:AC} we prove the following estimates for the inhomogeneous terms \(f,\tilde f\).

\lemma{InhomogeneousEst}{
We have the estimates
\begin{gather}
\|\Phi f\|_{\K}\lesssim \epsilon\label{est:InhomogeneousL2},\\
\|\Phi\tilde f\|_{\tilde \K}\lesssim\epsilon\label{est:InhomogeneousH01}.
\end{gather}
}
\bpf
We estimate each component of the \(Z,\tilde Z\) norms separately.

\textit{A. Estimating \(\|\Phi f\|_{L^\infty_TL^2}\).}
We start by calculating
\begin{align*}
f &= t^{-\frac{1}{3}}Q_w\W_t{+t^{-1}RQ_w\partial_x\W} + 6\sigma t^{-1}Q^2Q_w\partial_x\W + t^{-\frac{2}{3}}Q_{wy}\partial_x^2\W\\
&\quad + \tfrac{1}{3}t^{-\frac{1}{3}}Q_w\partial_x^3\W + t^{-\frac{2}{3}}Q_{wwy}(\partial_x\W)^2 + t^{-\frac{1}{3}}Q_{ww}\partial_x\W\partial_x^2\W + \tfrac{1}{3}t^{-\frac{1}{3}}Q_{www}(\partial_x\W)^3,
\end{align*}
where \(R(y)=\frac23y-\frac{\zeta(y)}{3\zeta'(y)}\) vanishes for \(|y|\geq1\) and we have used that
{
\[
\partial_t\W(t,t^{-\frac{1}{3}}\zeta(t^{-\frac13}x))=\left(t^{-\frac23}R(t^{-\frac13}x)-t^{-1}x\right)\partial_x(\W(t,t^{-\frac{1}{3}}\zeta(t^{-\frac13}x)))+\W_t(t,t^{-\frac{1}{3}}\zeta(t^{-\frac13}x)).
\]
}
We claim that
\est{Inhomogeneousf}{
\|f\|_{L^2}\lesssim \epsilon t^{-\frac{4+\delta}{3}},
}
and hence using \eqref{est:KPV}
\[
\|\Phi f\|_{L^\infty_TL^2}\lesssim\|f\|_{L^1([T,\infty);L^2)}\lesssim \epsilon T^{-\frac{1+\delta}{3}}.
\]

To prove \eqref{est:Inhomogeneousf} we estimate each term using Lemma \ref{lem:WEsts} for \(\W\) and Lemma \ref{lem:PLIIEsts} for \(Q\). For the first term we place \(\W_t\) into \(L^2\) and the remaining terms into \(L^\infty\) to get
\[
\|t^{-\frac13}Q_w\W_t\|_{L^2}\lesssim \|t^{-\frac13}\langle t^{-\frac13}x\rangle^{-\frac14 + \frac\delta4}\W_t\|_{L^2}\lesssim \epsilon t^{-\frac{4+\delta}{3}}.
\]
For the second term we use that \(R\) is a bounded function, supported in the region \(|y|\lesssim1\) to obtain
\[
\|t^{-1}RQ_w\partial_x\W\|_{L^2}\lesssim \|R\|_{L^2}\|t^{-1}\langle t^{-\frac13}x\rangle^{-\frac14 + \frac\delta4}\partial_x\W\|_{L^\infty}\lesssim\epsilon t^{-\frac{4+\delta}{3}}.
\]
For the third term we estimate in a similar manner to get
\[
\|t^{-1}Q^2Q_w\partial_x\W\|_{L^2}\lesssim \|t^{-1}\langle t^{-\frac13}x\rangle^{-\frac34+\frac\delta4}\|_{L^2}\|\W\|_{L^\infty}^2\|\partial_x\W\|_{L^\infty}\lesssim\epsilon^3 t^{-\frac{4+\delta}{3}}.
\]
For the remaining terms we place one \(\W\) into \(L^2\) and the rest into \(L^\infty\) to obtain
\begin{gather*}
\|t^{-\frac23}Q_{wy}\partial_x^2\W\|_{L^2}\lesssim \|t^{-\frac23}\langle t^{-\frac13}x\rangle^{\frac14+\frac\delta4}\partial_x^2\W\|_{L^2}\lesssim\epsilon t^{-\frac{4+\delta}{3}},\\
\|t^{-\frac13}Q_w\partial_x^3\W\|_{L^2}\lesssim \|t^{-\frac13}\langle t^{-\frac13}x\rangle^{-\frac14+\frac\delta4}\partial_x^3\W\|_{L^2}\lesssim\epsilon t^{-\frac{4+\delta}{3}},\\
\|t^{-\frac23}Q_{wwy}(\partial_x\W)^2\|_{L^2}\lesssim\|\W\|_{L^\infty}\|t^{-\frac23}\langle t^{-\frac13} x\rangle^{\frac14+\frac\delta4}\partial_x\W\|_{L^\infty}\|\partial_x\W\|_{L^2}\lesssim \epsilon^3 t^{-\frac32-\frac\delta3},\\
\|t^{-\frac13}Q_{ww}\partial_x\W\partial_x^2\W\|_{L^2}\lesssim\|\W\|_{L^\infty}\|\partial_x\W\|_{L^\infty}\|t^{-\frac13}\langle t^{-\frac13}x\rangle^{-\frac14+\frac\delta4}\partial_x^2\W\|_{L^2}\lesssim\epsilon^3t^{-\frac32-\frac\delta3},\\
\|t^{-\frac13}Q_{www}(\partial_x\W)^3\|_{L^2}\lesssim\|\partial_x\W\|_{L^\infty}^2\|t^{-\frac13}\langle t^{-\frac13}x\rangle^{-\frac14+\frac\delta4}\partial_x\W\|_{L^2}\lesssim\epsilon^3t^{-\frac{5+\delta}{3}}.
\end{gather*}
Combining these estimates we obtain \eqref{est:Inhomogeneousf}.

\textit{B. Estimating \(\|\Phi f_x\|_{L^4_xL^\infty_T}\).}
Computing \(f_x\) and estimating each term as before using Lemma \ref{lem:WEsts} for \(\W\) and Lemma \ref{lem:PLIIEsts} for \(Q\) we obtain
\est{fxNotEnough}{
\|f_x\|_{L^2}\lesssim \epsilon t^{-1-\frac\delta3},
}
whenever \(0<\epsilon\ll1\) is sufficiently small. Interpolating with the estimate \eqref{est:Inhomogeneousf} we obtain
\[
\||D|^{\frac14}f\|_{L^2}\lesssim \epsilon t^{-\frac54-\frac\delta3},
\]
and hence from the estimate \eqref{est:KPV} we have
\[
\|\Phi f\|_{L^4_xL^\infty_T}\lesssim \||D|^{\frac14}f\|_{L^1([T,\infty);L^2)}\lesssim \epsilon T^{-\frac14-\frac\delta3}.
\]

\textit{C. Estimating \(\|\Phi f_x\|_{L^\infty_TL^2_x}\).} Na\"ively proceeding as above, we have the estimate
\[
\|\Phi f_x\|_{L^\infty_TL^2_x}\lesssim \|f_x\|_{L^1([T,\infty);L^2)}\lesssim \delta^{-1}\epsilon T^{-\frac\delta3},
\]
which is insufficient to prove \eqref{est:InhomogeneousL2} as \(\delta\sim \epsilon^2\). Instead we first decompose
\(
f = g + b,
\)
into a good part \(g\) and a bad part \(b\) where,
\begin{align*}
g &= t^{-1}RQ_w\partial_x\W + 6\sigma t^{-1}Q^2Q_w\partial_x\W + \tfrac{1}{3}t^{-\frac{1}{3}}Q_w\partial_x^3\W + t^{-\frac{1}{3}}Q_{ww}\partial_x\W\partial_x^2\W\\
&\quad + t^{-\frac{2}{3}}Q_{wwy}(\partial_x\W)^2 + \tfrac{1}{3}t^{-\frac{1}{3}}Q_{www}(\partial_x\W)^3,\\
b &= t^{-\frac{1}{3}}Q_w\W_t  + t^{-\frac{2}{3}}Q_{wy}\partial_x^2\W.
\end{align*}
For the good part we may estimate as before using Lemmas \ref{lem:WEsts} and \ref{lem:PLIIEsts} to obtain the improved bound
\[
\|g_x\|_{L^2}\lesssim \epsilon t^{-\frac76 - \frac\delta3},
\]
whenever \(0<\epsilon\ll1\) is sufficiently small. As a consequence we have the estimate
\[
\|\Phi g_x\|_{L^\infty_TL^2}\lesssim\|g_x\|_{L^1([T,\infty);L^2)}\lesssim\epsilon T^{-\frac16-\frac\delta3}.
\]

For the bad part we will use that frequency localization of \(S(-t)\uapp\) will correspond to spatial localization of \(\W\). As a consequence we will aim to show that
\est{bxGoal}{
\|S(-t)b_x\|_{l^2L^1([T,\infty);L^2)} \lesssim\epsilon T^{-\frac\delta3},
}
where the \(l^2\)-summation is with respect to dyadic regions in \(x\)-frequency. We may then use the embedding \eqref{est:UpVpEmbed} and that the \(V^2\) norm commutes with the \(l^2\)-summation to obtain
\[
\|\Phi b_x\|_{L^\infty_TL^2}\lesssim \|\Phi b_x\|_{V^2([T,\infty);L^2)}\lesssim \|\Phi b_x\|_{l^2V^2([T,\infty);L^2)}\lesssim \|S(-t)b_x\|_{l^2L^1([T,\infty);L^2)}\lesssim \epsilon T^{-\frac\delta3}.
\]

To prove \eqref{est:bxGoal}, we first note that we have improved bound for low frequencies,
\[
\|P_{\leq T^{\frac13}}\partial_xS(-t)b\|_{L^2}\lesssim T^{\frac13}\|b\|_{L^2}\lesssim\epsilon T^{\frac13}t^{-\frac{4+\delta}{3}}.
\]
Integrating we obtain
\[
\|P_{\leq T^{\frac13}}S(-t)b_x\|_{L^1([T,\infty);L^2)}\lesssim \epsilon T^{-\frac\delta3}.
\]

For a dyadic frequency \(M>T^{\frac13}\) and \(t\geq T\) we expect that \(P_Mb\) will be localized in the spatial region \(\{|x|\sim tM^2\}\). Taking \(\chi_M\) as in \S\ref{sect:InitBnds} and applying the elliptic estimate \eqref{est:Elliptic}, we have
\[
\|P_Mb\|_{L^2}\lesssim \|P_M(\chi_Mb)\|_{L^2} + t^{-1}M^{-2}\|Lb\|_{L^2} + t^{-1}M^{-3}\|b\|_{L^2}.
\]
Again applying the estimates of Lemmas \ref{lem:WEsts} and \ref{lem:PLIIEsts}, we obtain the bound
\[
\|Lb\|_{L^2}\lesssim \epsilon t^{-1-\frac\delta3}.
\]
As a consequence we have
\[
\|P_MS(-t)b_x\|_{L^2} \lesssim M\|P_M(\chi_Mb)\|_{L^2} + \epsilon t^{-2-\frac\delta3}M^{-1},
\]
and hence
\[
\int_T^\infty \|P_MS(-t)b_x\|_{L^2}\,dt\lesssim \int_T^\infty M\|P_M(\chi_Mb)\|_{L^2}\,dt + \epsilon T^{-1-\frac\delta3}M^{-1}.
\]
From the proof of Lemma \ref{lem:WEsts} we have the slightly refined estimates
\begin{align}
\|\chi_M b\|_{L^2} &\lesssim t^{-\frac32 - \frac\delta3}M^{-\frac12}\|\chi_M\langle D\rangle^{1+\delta}W\|_{L^2} + \epsilon t^{-\frac32-\frac\delta3}M^{-\frac32},\label{est:chimb}\\
\|\partial_x^2(\chi_M b)\|_{L^2} &\lesssim t^{-\frac\delta3}M\|\chi_M\langle D\rangle^\delta W\|_{L^2} + \epsilon t^{-\frac12-\frac\delta3}M^{-\frac12},\label{est:ddchimb}
\end{align}
where we have used the fact that \(\W\) only depends on \(W\)-frequencies \(\leq t\) in the second estimate. For large times we use \eqref{est:chimb} to obtain
\begin{align*}
\int_{\max\{M,T\}}^\infty M\|P_M(\chi_Mb)\|_{L^2}\,dt &\lesssim \int_{\max\{M,T\}}^\infty \left(t^{-\frac32 - \frac\delta3}M^{\frac12}\|\chi_M\langle D\rangle^{1+\delta}W\|_{L^2} + \epsilon t^{-\frac32-\frac\delta3}M^{-\frac12}\right)\,dt\\
&\lesssim T^{-\frac\delta3}\|\chi_M\langle D\rangle^{1+\delta}W\|_{L^2} + \epsilon T^{-\frac12-\frac\delta3}M^{-\frac12}.
\end{align*}
If \(M>T\) we use the estimate \eqref{est:ddchimb} to obtain
\begin{align*}
\int_T^M M\|P_MS(-t)(\chi_Mb)\|_{L^2}\,dt &\lesssim \int_T^M \left(t^{-\frac\delta3}\|\chi_M\langle D\rangle^\delta W\|_{L^2} + \epsilon t^{-\frac12-\frac\delta3}M^{-\frac32}\right)\,dt\\
&\lesssim T^{-\frac\delta3}M\|\chi_M\langle D\rangle^\delta W\|_{L^2} + \epsilon T^{-\frac12-\frac\delta3}M^{-\frac12}.
\end{align*}
Summing these we have
\[
\int_T^\infty \|P_MS(-t)b_x\|_{L^2}\,dt\lesssim T^{-\frac\delta3}\|\chi_M\langle D\rangle^{1+\delta}W\|_{L^2} + T^{-\frac\delta3}M\|\chi_M\langle D\rangle^\delta W\|_{L^2} + \epsilon T^{-\frac12-\frac\delta3}M^{-\frac12}.
\]
Summing dyadyically over \(M>T^{\frac13}\) we obtain the estimate \eqref{est:bxGoal}.

\textit{D. Estimating \(\|\Phi\tilde f\|_{L^\infty_TL^2}\).}
From the estimate \eqref{est:uappBndsInf} for \(\uapp\) and \eqref{est:Inhomogeneousf} for \(f\) we have
\[
\|t\uapp^2f\|_{L^2}\lesssim t\|\uapp\|_{L^\infty}^2\|f\|_{L^2}\lesssim\epsilon^3 t^{-1-\frac\delta3},
\]
and hence
\[
\|t\uapp^2f\|_{L^1([T,\infty);L^2)}\lesssim\epsilon T^{-\frac\delta3}.
\]

To estimate \(Lf\) we decompose \(f = g + b + r\) into a good part \(g\), a bad part \(b\) and a ``cubic'' part \(r\) defining
\begin{align*}
g &= t^{-\frac{2}{3}}Q_{wwy}(\partial_x\W)^2 + t^{-\frac{1}{3}}Q_{ww}\partial_x\W\partial_x^2\W + \tfrac{1}{3}t^{-\frac{1}{3}}Q_{www}(\partial_x\W)^3,\\
b &= t^{-\frac{1}{3}}Q_w\W_t +t^{-1}RQ_w\partial_x\W + t^{-\frac{2}{3}}Q_{wy}\partial_x^2\W + \tfrac{1}{3}t^{-\frac{1}{3}}Q_w\partial_x^3\W,\\
r &= 6\sigma t^{-1}Q^2Q_w\partial_x\W.
\end{align*}

Estimating as before using Lemmas \ref{lem:WEsts} and \ref{lem:PLIIEsts} we obtain the improved bound
\[
\|Lg\|_{L^2}\lesssim\epsilon^3 t^{-1-\frac\delta3},
\]
which we may integrate in time to obtain
\[
\|Lg\|_{L^1([T,\infty);L^2)}\lesssim\epsilon T^{-\frac\delta3}.
\]

For the bad piece we will use that spatial localization of \(S(-t)\uapp\) will correspond to frequency localization of \(\W\). In this case we aim to show that
\est{bxGoal2}{
\|S(-t)Lb\|_{l^2L^1([T,\infty);L^2)} \lesssim\epsilon T^{-\frac\delta3},
}
where the \(l^2\)-summation is now with respect to dyadic spatial regions in \(x\). As before we may commute the \(V^2\)-norm with the \(l^2\)-summation to obtain
\[
\|\Phi Lb\|_{L^\infty_TL^2}\lesssim\|S(-t)Lb\|_{l^2L^1([T,\infty);L^2)} \lesssim\epsilon T^{-\frac\delta3}.
\]

To prove \eqref{est:bxGoal2} we first observe that we have an improved estimate in the self-similar region,
\[
\|S(-t)Lb\|_{L^2(|x|\leq T^{\frac13})}\lesssim \|xS(-t)b\|_{L^2(|x|\leq T^{\frac13})}\lesssim T^{\frac13}\|b\|_{L^2}\lesssim \epsilon T^{\frac13}t^{-\frac{4+\delta}{3}}.
\]
Integrating we obtain
\[
\|\chi_{\{|x|\leq T^{\frac13}\}}S(-t)Lb\|_{L^1([T,\infty);L^2)}\lesssim \epsilon T^{-\frac\delta3}.
\]

From the proof of Lemma \ref{lem:WEsts} we once again obtain refined estimates
\begin{align}
\|b\|_{L^2} &\lesssim t^{-\frac32-\frac\delta3}\|W_{\leq t^{\frac13}}\|_{H^{\frac32+\delta}} + t^{-\frac76-\frac\delta3}\|W_{>t^{\frac13}}\|_{H^{\frac12+\delta}},\label{est:BadHelp1}\\
\|L^2b\|_{L^2} &\lesssim t^{-\frac56-\frac\delta3}\|W_{\leq t^{\frac13}}\|_{H^{\frac32+\delta}} + t^{-\frac12-\frac\delta3}\|W_{>t^{\frac13}}\|_{H^{\frac12+\delta}}\label{est:BadHelp2}.
\end{align}
Applying the Cauchy-Schwarz inequality with \eqref{est:BadHelp1} we obtain
\[
\left(\int_T^{M^3} \|\chi_{\{|x|\sim M\}}xS(-t)b\|_{L^2}\,dt\right)^2 \lesssim T^{-\frac{2\delta}{3}}\int_{M^3}^\infty M^{\frac12}t^{-\frac32}\left(\|W_{\leq t^{\frac13}}\|_{H^{\frac32+\delta}}^2 + t^{\frac23}\|W_{>t^{\frac13}}\|_{H^{\frac12+\delta}}^2\right)\,dt.
\]
Similarly, applying the Cauchy-Schwarz inequality with \eqref{est:BadHelp2} we obtain
\[
\left(\int_{M^3}^\infty \|\chi_{\{|x|\sim M\}}xS(-t)b\|_{L^2}\,dt\right)^2 \lesssim T^{-\frac{2\delta}{3}}\int_T^{M^3} M^{-\frac12}t^{-\frac76}\left(\|W_{\leq t^{\frac13}}\|_{H^{\frac32+\delta}}^2 + t^{\frac23}\|W_{>t^{\frac13}}\|_{H^{\frac12+\delta}}^2\right)\,dt.
\]
Summing these we have the estimate
\[
\|\chi_{\{|x|\sim M\}}xS(-t)b\|_{L^1([T,\infty);L^2)}^2\lesssim T^{-\frac{2\delta}{3}}\int_T^\infty \min\{M^{\frac12}t^{-\frac32},M^{-\frac12}t^{-\frac76}\}\left(\|W_{\leq t^{\frac13}}\|_{H^{\frac32+\delta}}^2 + t^{\frac23}\|W_{>t^{\frac13}}\|_{H^{\frac12+\delta}}^2\right)\,dt.
\]
We may then sum over dyadic \(M>T^{\frac13}\) to obtain
\[
\|xS(-t)b\|_{l^2L^1([T,\infty);L^2)}^2\lesssim T^{-\frac{2\delta}{3}}\int_T^\infty t^{-\frac43}\|W_{\leq t^{\frac13}}\|_{H^{\frac32+\delta}}^2 + t^{-\frac23}\|W_{>t^{\frac13}}\|_{H^{\frac12+\delta}}^2\,dt.
\]

Finally we decompose \(W\) by dyadic frequencies to calculate the integral
\begin{align*}
\int_T^\infty t^{-\frac43}\|W_{\leq t^{\frac13}}\|_{H^{\frac32+\delta}}^2 + t^{-\frac23}\|W_{>t^{\frac13}}\|_{H^{\frac12+\delta}}^2\,dt & \lesssim \|W_{\leq 1}\|_{L^2}^2 + \int_T^\infty \sum\limits_{N\geq1}\min\{N^{3+2\delta}t^{-\frac43},N^{1+2\delta}t^{-\frac23}\}\|W_N\|_{L^2}^2\,dt\\
&\lesssim \|W_{\leq 1}\|_{L^2}^2 + \sum\limits_{N\geq1}\int_T^\infty \min\{N^{3+2\delta}t^{-\frac43},N^{1+2\delta}t^{-\frac23}\}\|W_N\|_{L^2}^2\,dt\\
&\lesssim \|W_{\leq 1}\|_{L^2}^2 + \sum\limits_{N\geq1} N^{2(1+\delta)}\|W_N\|_{L^2}^2,
\end{align*}
which completes the proof of \eqref{est:bxGoal2}.

While the cubic piece \(r\) should be well-behaved as it is cubic in \(\W\), using Lemmas \ref{lem:WEsts} and \ref{lem:PLIIEsts} as before we can only obtain the estimate
\[
\|Lr + 12\sigma t^{-\frac23}(Q_y^2Q_w + QQ_{yy}Q_w + QQ_yQ_{wy})\partial_x\W\|_{L^2}\lesssim\epsilon^3 t^{-1-\frac\delta3}.
\]
The difficulty with the remaining terms in \(r\) is that they do not decay sufficiently in \(y\) to allow us to estimate \(\partial_x\W\) in \(L^\infty\). Fortunately however, the troublesome terms are non-resonant and may be removed by what is is essentially a normal form. Using that
\[
\partial_y^{k+2}\partial_w^jQ = y\partial_y^k\partial_w^jQ + \textrm{lower order terms}
\]
we define
\[
q(t,x) = \frac{9\sigma t^{\frac13}}{5\zeta(t^{-\frac13}x)^2}\left(2QQ_yQ_w + Q^2Q_{wy}\right)\partial_x\W.
\]
Using Lemmas \ref{lem:WEsts} and \ref{lem:PLIIEsts} we have
\[
\|q\|_{L^2}\lesssim\epsilon^3 t^{-\delta}.
\]
We calculate
\[
(\partial_t + \tfrac13\partial_x^3)q = - 12\sigma t^{-\frac23}(Q_y^2Q_w + QQ_{yy}Q_w + QQ_yQ_{wy})\partial_x\W + \err,
\]
where the error term \(\err\) may be estimated using Lemmas \ref{lem:WEsts} and \ref{lem:PLIIEsts} to obtain
\[
\|\err\|_{L^2}\lesssim\epsilon^3 t^{-1-\delta}.
\]
As a consequence we obtain the estimate
\[
\|\Phi Lr\|_{L^\infty_TL^2}\lesssim \|q\|_{L^\infty_TL^2} + \int_T^\infty \epsilon^3 t^{-1-\delta}\,dt\lesssim\epsilon T^{-\delta},
\]
which completes the proof of \eqref{est:InhomogeneousH01}.

\epf

\medskip


\section*{Acknowledgements}
The author would like to thank his advisor Daniel Tataru for suggesting the problem and several key suggestions for the proof. He would like to thank Herbert Koch for several useful comments and Mihaela Ifrim for a number of helpful discussions. He would also like to thank the reviewer for many helpful comments on the manuscript.

Part of this research was completed while the author attended the Hausdorff trimester program on ``Harmonic Analysis and Partial Differential Equations'' at the Hausdorff Research Institute for Mathematics. The author would like to thank the Hausdorff Institute for their hospitality and support during the trimester.


\begin{appendix}

\section{Short-range perturbations}\label{app:perturbations}

In this appendix we briefly outline some modifications to Theorems \ref{thrm:Main} and \ref{thrm:AC} in the case of short-range perturbations \eqref{eq:PerturbedmKdV}.

When \(p\in[\tfrac{7}{2},\infty)\) the results are essentially unchanged. For \(p\in(3,\tfrac{7}{2})\) the energy estimate \eqref{est:JNRG} fails for \(\delta\) defined as in \eqref{defn:delta}. This is due to the fact that
\(v=\Lambda u\) satisfies the equation
\[
\pde{
v_t+\tfrac{1}{3}v_{xxx}=(3\sigma u^2+F'(u))v_x+3F(u)-F'(u)u,
}{
v(0)=xu_0.
}
\]
So, if \(\epsilon>0\) is sufficiently small, estimating \(\|\Lambda v\|_{L^2}\) as in \eqref{est:JNRG} we have to re-define
\defn{newDelta}{
\delta=\frac{7-2p}{6}\in(0,\tfrac{1}{6}).
}
In particular, for sufficiently small \(\epsilon>0\), the loss of regularity in \eqref{est:AEst} is controlled by \(p\) rather than \(\epsilon\), giving us the revised estimate
\[
\|W\|_{H^{1-C\delta,1}\cap L^\infty}\lesssim\epsilon.
\]

\end{appendix}

\bibliographystyle{abbrv}
\bibliography{mKdVSubmission}

\end{document}